\documentclass[11pt]{article}

\usepackage{amsbsy,amstext,amsopn,amssymb}
\usepackage{latexsym}
\usepackage{amscd}
\usepackage{amsmath}
\usepackage{amsfonts}
\usepackage{theorem}
\usepackage{graphicx}
\usepackage{color}

{\theorembodyfont{\slshape}
\newtheorem{theorem}{Theorem}
\newtheorem{proposition}{Proposition}
\newtheorem{lemma}{Lemma}

} {\theorembodyfont{\rmfamily}
\newtheorem{definition}{Definition}

\newtheorem{remark}{Remark}

}

\newcommand{\matr}[1]{\begin{bmatrix} #1 \end{bmatrix}}    

\newcommand{\gscmx}{{\chi}}

\def\proof{{\noindent\sc Proof. \quad}}
\newcommand\proofof[1]{{\noindent\sc Proof of #1. \quad}}
\def\eproof{{\mbox{}\hfill\qed}\medskip}
\newcommand\qed{{\unskip\nobreak\hfil\penalty50\hskip2em\vadjust{}
\nobreak\hfil$\Box$\parfillskip=0pt\finalhyphendemerits=0\par}}

\def\CACM{Communications of the ACM}

\def\SIOPT{SIAM Journal on Optimization}

\def\JoC{J. of Complexity}
\def\MP{Math. Program.}
\begin{document}

\bibliographystyle{plain}

\makeatletter

\renewcommand{\bar}{\overline}

\def\transp{^{\rm T}}

\newcommand{\x}{\times}
\newcommand{\<}{\langle}
\renewcommand{\>}{\rangle}
\newcommand{\into}{\hookrightarrow}

\renewcommand{\a}{\alpha}
\renewcommand{\b}{\beta}
\renewcommand{\d}{\delta}
\newcommand{\D}{\Delta}
\renewcommand{\th}{\theta}
\newcommand{\g}{\gamma}
\renewcommand{\l}{\lambda}
\renewcommand{\t}{\tau}
\renewcommand{\O}{\Omega}
\newcommand{\z}{\zeta}

\newcommand{\p}{\partial}
\renewcommand{\hat}{\widehat}
\newcommand{\ubar}{\underline}
\renewcommand{\tilde}{\widetilde}

\newcommand{\constCQR}{\mathbf{c}_1}
\newcommand{\constCL}{\mathbf{c}_2}
\newcommand{\constCGolub}{\mathbf{c}_3}


\font\eightrm=cmr8 \font\ninerm=cmr9


\makeatletter
\def\R{\mathchoice{{\setbox0=\hbox{\rm I}\copy0\kern-0.55\wd0\hbox{\rm R}}}{%
  {\setbox0=\hbox{\rm I}\copy0\kern-0.55\wd0\hbox{\rm R}}}{%
  {\setbox0=\hbox{$\m@th\scriptstyle\rm I$}\copy0\kern-0.5\wd0%
  \scriptstyle\rm R}}{%
  {\setbox0=\hbox{$\m@th\scriptscriptstyle\rm I$}\copy0\kern-0.6\wd0%
  \scriptscriptstyle\rm R}}}
\def\N{\mathchoice{{\setbox0=\hbox{\rm I}\copy0\kern-0.55\wd0\rm N}}{%
  {\setbox0=\hbox{\rm I}\copy0\kern-0.55\wd0\rm N}}{%
  {\setbox0=\hbox{$\m@th\scriptstyle\rm I$}\copy0\kern-0.5\wd0%
  \scriptstyle\rm N}}{%
  {\setbox0=\hbox{$\m@th\scriptscriptstyle\rm I$}\copy0\kern-0.6\wd0%
  \scriptscriptstyle\rm N}}}
\def\F{\mathchoice{{\setbox0=\hbox{\rm I}\copy0\kern-0.55\wd0\hbox{\rm F}}}{%
  {\setbox0=\hbox{\rm I}\copy0\kern-0.55\wd0\hbox{\rm F}}}{%
  {\setbox0=\hbox{$\m@th\scriptstyle\rm I$}\copy0\kern-0.5\wd0%
  \scriptstyle\rm F}}{%
  {\setbox0=\hbox{$\m@th\scriptscriptstyle\rm I$}\copy0\kern-0.6\wd0%
  \scriptscriptstyle\rm F}}}
\def\Oh{{\cal O}}
\def\Error{{\hbox{\tt Error}}}
\def\fl{{\hbox{\tt fl}}}
\def\bfc{{\bf c}}
\def\bn{{\bf n}}
\def\bm{{\bf m}}
\def\AA{{\cal A}}
\def\BB{{\cal B}}
\def\KK{{\cal{K}}}
\def\LL{{\cal L}}
\def\NN{{\cal N}}
\def\rou{\mbox{\sf round}}
\renewcommand{\int}{\mbox{\sf int}}
\renewcommand{\det}{\mbox{\sf det\,}}
\def\vx{\vec{x}}
\def\vy{\vec{y}}
\def\vs{\vec{s}}
\def\Lg{{\rm Lg}\,}
\def\Id{{\rm I}}
\def\bfr{{\bf r}}
\def\um{{\scriptscriptstyle 1/2}}
\newcommand{\refer}{\eqref}
\newcommand{\ip}[2]{\langle #1 , #2 \rangle}    
\def\ll{{[\kern-1.6pt [}}
\def\rr{{]\kern-1.4pt ]}}

\newcommand{\ch}{\textcolor{red}}

\begin{title}
{{\bf  Solving second-order conic systems with variable precision}}
\end{title}
\author{
Felipe Cucker\thanks{Partially supported by GRF grant
CityU 1008/08}\\
Department of Mathematics\\
City University of Hong Kong\\
83 Tat Chee Avenue, Kowloon\\
HONG KONG\\
e-mail:
{\tt macucker@math.cityu.edu.hk}
\and
Javier Pe\~na\thanks{Supported by NSF grant
CCF-0830533}\\
Tepper School of Business\\
Carnegie Mellon University\\
5000 Forbes Avenue, Pittsburgh, PA 15213-3890\\
USA\\
e-mail: {\tt jfp@andrew.cmu.edu}
\and
Vera Roshchina\\
Collaborative Research Network\\
University of Ballarat\\
Mount Helen Campus\\
PO Box 663, CRN F Building\\
VIC 3353 Australia\\
e-mail:
{\tt vroshchina@ballarat.edu.au}}

\date{}

\makeatletter \maketitle \makeatother

\begin{quote}
{\small {\bf Abstract.} We describe and analyze an interior-point
method to decide feasibility problems of second-order conic
systems. A main feature of our algorithm is that arithmetic
operations are performed with finite precision. Bounds for
both the number of arithmetic operations and the finest precision
required are exhibited.}
\end{quote}


\section{Introduction}\label{s:intro}

It is now widely accepted that the most efficient algorithms for
solving the general type of second-order conic problems are
interior-point methods (IPMs). IPMs infallibly demonstrate very
fast numerical convergence, by far outperforming their theoretical
estimates.

Second-order conic programming problems contain linear programming
problems as a special case, and at the same time can be embedded
into the class of semidefinite programming problems. It is,
however, not advisable to solve SOCP problems by semidefinite
programming methods (see \cite{AlizGold03}, \cite{SVBL:98}) as
IPMs that solve SOCP directly have a much better complexity (both
in theory and in practice). SOCP problems have lately received
considerable attention due to their many applications
\cite{SVBL:98}; they appear to be at the boundary of the problems
for which interior-point methods can solve large instances, a fact
that is linked to the implementation of commercial software for the
solution of second-order programs such as MOSEK\footnote{\tt
http://www.mosek.com/} or CPLEX\footnote{\tt
http://www.ilog.com/products/cplex/}.

We are interested in solving homogeneous second-order conic feasibility problems.
That is, given a second-order cone $K\subset \R^{n}$ and a matrix $A\in
\R^{m\times n}$\label{glossary:A}, decide which one of the primal-dual pair of
problems
$$
  \begin{array}{ll}
      Ax=0,\\
      x\succeq_K 0,
   \end{array}
   \quad({\rm P})
   \qquad\qquad
   \begin{array}{ll}
     A\transp y+s =0, \\
     s\succeq_K 0,
     \end{array}
     \quad({\rm D})
$$
is strictly feasible (i.e., the relevant conic constraint is
strict) and provide a solution to the feasible problem. It is
well-known that each of (P) and (D) above has a strict solution if
and only if the other one has no nonzero solutions.

Recall that a second-order cone is a direct product of a finite
number of Lorentz cones. The {\em Lorentz cone} $\LL_p\subset \R^{p+1}$\label{glossary:LL}
is defined to be
$$
    \LL_p = \{x\in \R^{p+1}\,|\, x_0\geq \|\bar x\|\},
$$
where for a vector $x\in \R^{p+1}$ indexed from $0$ to $p$ we
let $\bar x = (x_1, x_2, \dots, x_p)\in\R^p$. For our
primal-dual pair of problems (P)-(D) we take
$$
   K = \LL_{n_1}\times \LL_{n_2}\times\dots\times \LL_{n_r},\label{glossary:K}
$$
where $r$\label{glossary:rplain} is the number of  Lorentz cones comprising $K$, and
$\sum_{i=1}^{r}(n_i+1) = n$ with $n_i$ being positive integers for all $i$ from
1 to $r$.

We propose a finite-precision algorithm for solving the SOCP
feasibility problem and provide rigorous bounds for the finest
machine precision and the maximal number of iterations needed.
The proposed algorithm
is designed to work with variable precision, that is, the machine
precision can be re-adjusted along the way.

Our bounds depend on Renegar's condition
number~\cite{Renegar95b},\cite{Renegar99}, which is consistent with
similar bounds obtained for the polyhedral case in \cite{CP01}. Let
$\rho_P(A)$ and $\rho_D(A)$ be the {\em distance to infeasibility}
of (P) and (D) respectively defined by
$$
  \rho_P(A)=\inf\{\|\Delta A\|:(A+\Delta
  A)x=0,x\succ_{K}0 \mbox{ is infeasible}\}\label{glossary:rhoP}
$$
and
$$
  \rho_D(A)=\inf\{\|\Delta A\|:-(A+\Delta A)\transp y\succ_{K}0,
  y\in\R^m \mbox{ is infeasible}\}.\label{glossary:rhoD}
$$
Renegar's {\em condition number} $C(A)$ is defined as the reciprocal of the relative distance to ill-posedness of the pair
(P)--(D):
$$
   C(A):=\frac{\|A\|}{\max\{\rho_P(A),\rho_D(A)\}}.\label{glossary:CA}
$$
Although any equivalent matrix norm can be used to define $C(A)$, in
our analysis we choose to use the standard operator norm induced by
the Euclidean scalar product. We say that the problem is ill-posed if
both $\rho_P(A) =\rho_D(A) = 0$ and hence $C(A) = \infty$.

Our main result, Theorem~\ref{main_th}, shows
that there exists a finite precision interior-point
method which, with input a matrix $A\in \R^{m\times n}$ and a
second-order conic structure $K$ (consisting of $r$ Lorentz cones),
decides which one of the two systems (P) or (D) is feasible. We
estimate both the number of iterations of the algorithm and the
precision required as functions of the size of the matrix, the
number $r$ of Lorentz cones in $K$, and the condition number of the
problem. The finest required precision is
$$
   u = \frac{1}{\Oh\bigl((m+n)^{5/2} r^{11.5} C(A)^{7/2}\bigr)},\label{glossary:Oh}
$$
and the number of main interior-point iterations performed by
the algorithm is bounded by
$$
  \Oh(r^\um (\log r+\log C(A))).
$$

Strictly speaking, our algorithm solves both the {\em decision
problem} ---decide which one of the problems (P) and (D) is feasible---
and {\em the function problem} ---if either one of the problems (P) or
(D) is strictly feasible produce a (possibly approximate) solution for it.
For a precise version of our main result, the reader should check 
the statement of Theorem~\ref{main_th}.

Throughout the paper, we use standard notation wherever possible.
We index our variables according to the second-order conic
structure. That is, $x = (x_1, x_2, \dots , x_r)$, $s = (s_1, s_2,
\dots , s_r)$, where $x_i,s_i\in \R^{n_i+1}$ for all $i = 1,
\dots, r$. Throughout the paper we assume that
$\|A_i\|_F=1/\sqrt{r}$, and hence $\|A\|_F=1$. Note that this
assumption is trivial from a computational viewpoint; if $A_i \neq
0_{m\times (n_i+1)}$, it takes a few operations to reduce the matrix
to this form and it is easy to recover solutions of the original
system from those for the reduced one. The condition number of the
new matrix may change, however. But one can show as in
\cite[\S11.2]{CP01} that this change can not be large.

Finite precision analyses are pervasive in Numerical Linear Algebra;
they are much less common in optimization. While the effects of
finite precision when solving linear
programming problems had been early noticed
(e.g.\ \cite{Bartels,Clasen,Ogr,Robinson,Storoy,Wolfe})
there was no condition-based round-off
analysis even for linear programming problems until recently.
This was done for the feasibility problem for
polyhedral conic systems~\cite{CP01}, for the optimal
value of linear programs~\cite{Vera98}, and for the
computation of optimal basis and optimal solutions of linear
programs~\cite{ChC03}. To the best of our knowledge,
our work is the first such analysis for
nonlinear cones.

Our paper is organized as follows. In Section~\ref{s:ipm} we use a
relaxation scheme introduced by Pe\~na and Renegar~\cite{PeRe} and Vera et al. \cite{VeraEtAl2007} to reformulate the feasibility problem via an optimization one and recall the
basic ideas of interior-point methods. Then we relax the
standard results of IPM analysis to make room for computational
errors. We do not deal with finite-precision issues directly until
Section~\ref{s:fp}, where we describe our algorithm in detail and
estimate errors arising on every step of floating-point
computations. The last section is devoted to the proof of the main
result, and essentially fits the error estimates obtained in
Section~\ref{s:fp} into the gaps made for this purpose in
our extension of the IPM analysis done in Section~\ref{s:ipm}.

\section{Interior-point method for SOCP feasibility problem}
\label{s:ipm}

We use a relaxation scheme introduced by Pe\~na and Renegar in
\cite{PeRe} and later extended in \cite{VeraEtAl2007}. This relaxation scheme reformulates the feasibility problem (P)--(D), for the more general case when $K$ is a symmetric cone, as a pair of primal-dual optimization problems in higher dimension and solves this pair by a standard short-step interior-point method.  We next summarize the main ingredients of this approach.

It was shown in \cite{VeraEtAl2007} that the pair (P)--(D) is
equivalent to the following primal-dual pair of optimization problems
\begin{align}\label{eq:OP}
   \min\ &\vec{c}\transp\vec{x}\notag\\
   \mbox{s.t. }&\AA\vec{x}=\vec{b}\tag{{P'}}\\
   &\vec{x}\succeq_\KK 0\notag
 \end{align}
and
\begin{align}\label{eq:OD}
  \max\ &\vec{b}\transp\vec{y}\notag\\
  \mbox{s.t. }&\AA\transp\vec{y}+\vec{s}=\vec{c}\tag{{D'}}\\
  &\vec{s}\succeq_\KK 0\notag,
\end{align}
where $\KK = K\times \LL_n\times \LL_m$\label{glossary:KK}, $\vx = (x,t,x',\t,x'')$\label{glossary:vx}, 
$\vs = (s,t_s,s',\t_s,s'')\in \R^{n+1+n+1+m}$\label{glossary:vs}, $\vy = (y,y',-\eta)\in \R^{m+n+1}$\label{glossary:vy},
$$
   \AA:=\left[\begin{array}{ccccc} {A}&0&0&0&\Id_m\\
   -\Id_{n}&0&\Id_{n}&0&0\\
   0&1&0&0&0\end{array}\right],\quad  \vec{b}=\left[\begin{array}{c} 0\\
   0\\
  1\end{array}\right],\quad \vec{c}=\left[\begin{array}{c}0\\
  0\\
  0\\
  1\\
  0\end{array}\right].\label{glossary:AA}
$$
This equivalence should be understood in the following sense: If
$\rho(A) >0$ then a primal-dual interior-point method applied to the
pair~\eqref{eq:OP}--\eqref{eq:OD} yields a strict solution to
whichever of (P) or (D) is strictly feasible. In particular, since the
optimal value of the pair~\eqref{eq:OP}--\eqref{eq:OD} is zero, a
corresponding strict solution to the original feasibility problem can be
straightforwardly recovered from the first entries of the extended
variables ($x$ in the case of (P) and $y$ and $s$ in the case of (D)).

\ch{Note that from $\|A\|_F = 1$ we get
\begin{equation}\label{eq:normAA}
\|\AA\|  =
 \max_{\|\vec x\|\leq 1}\sqrt{\|Ax-x''\|^2+\|x'-x\|^2+t^2}
  \leq \max_{\|\vec x\|\leq
1}\sqrt{2\|Ax\|^2+2\|\vec x\|^2} \leq 2.
\end{equation}
}

In the sequel, to simplify notation, we will denote $\bm:=m+n+1$\label{glossary:bm}
and $\bn:=2n+m+2$\label{glossary:bn} so that $\AA\in\R^{\bm\times\bn}$,
$\vec{b}\in\R^\bm$ and $\vec{c}\in\R^\bn$. We also let
$\bfr:= r+2$\label{glossary:bfr}.

The pair \eqref{eq:OP}--\eqref{eq:OD} can be solved
via a primal-dual interior-point algorithm. We refer the reader to
\cite{Renegar99} and the references therein for a detailed exposition
of the theory of IPMs.  We next recall the concepts and results from
this basic theory that will be used in the paper.  Consider the
following self-scaled barrier function for the cone $\KK$:
\begin{equation}\label{eq:f}
f(\vec x) = -\left(\sum_{i=1}^{r} \ln (x_{0i}^2-\|\bar
{x_i}\|^2)+\ln (t^2-\|x'\|^2)+\ln (\t^2-\|x''\|^2)\right).
\end{equation}
Let $g(\vx) = \nabla f(\vx)$\label{glossary:grad} and $H(\vx) = \nabla^2 f(\vx)$\label{glossary:hess} denote
respectively the gradient and the Hessian of $f$. Let $e \in \KK$
denote the unique point such that $H(e) = I$, that is, $e =
(e_1,\dots,e_\bfr)$\label{glossary:e} where $e_{i0} = 1$ and $\bar e_i = 0$ for
$i=1,\dots,\bfr$.

Sometimes we will also need to work with the self-scaled barrier
function for the cone $K$:
\[
\bar f(x) = -\sum_{i=1}^{r} \ln (x_{0i}^2-\|\bar
{x_i}\|^2),\label{glossary:barf}
\]
and will let $\bar g(x) = \nabla \bar f(x)$\label{glossary:barg} and $\bar H(x) =
\nabla^2 \bar f(x)$\label{glossary:barH}.  Notice that $H(\vx) = \matr{\bar H(x) & 0
\\ 0 & \tilde H(\vx)}$, where $\tilde H(\vx)$ denotes the
Hessian of the function $-\ln (t^2-\|x'\|^2)-\ln
(\t^2-\|x''\|^2)$.

Given $\vx \in \KK$, the {\em local norm} $\|\cdot\|_{\vx}$ in
$\R^\bn$ is defined as
\[
 \|u\|_{\vx}^2 := u \transp H(\vx) u.
\]
Likewise for $x \in K$.

The {\em central path} of \eqref{eq:OP}--\eqref{eq:OD} is the set of
solutions $\{(\vec x(\mu), \vec y(\mu), \vec s(\mu))\in\int(\KK)\times\R^{\bm }\times\int(\KK): \mu > 0\}$
to the system of equations
\begin{eqnarray}\label{eq:central}
 \AA\vec{x}&=&\vec{b}\notag\\ \AA\transp\vec{y}+\vec{s}&=&\vec{c}\\
 \vec{s}+\mu g(\vec{x})&=&0.\notag
\end{eqnarray}
Given $z=(\vx,\vy,\vs)\in \int(\KK)\times\R^{\bm }\times\int(\KK)$\label{glossary:zplain} define
$$
  \mu(z) := \frac{\vec{x}\transp \vec{s}}{2\bfr}.\label{glossary:muofz}
$$
Note that if $z$ belongs to the central path for a certain value of
$\mu$ then $\mu(z)=\mu$. We may sometimes write $\mu$ for $\mu(z)$
when $z$ is clear from the context.

The basic idea of a path-following interior-point method is to
generate a sequence of points on a suitable neighborhood of the
central path that converges to optimality.  The suitable neighborhood
is the following.

\begin{definition}
Given \label{glossary:b} $\b\in(0,1/15)$, the {\em central neighborhood}
$\NN_\b$\label{glossary:NNb} is defined as the set of points $z=(\vec{x},\vec{y},\vec{s})\in \int(\KK)\times\R^{\bm }\times\int(\KK)$,
such that the following constraints hold:
$$
  \begin{array}{rl} \AA\vec{x}&=\;\vec{b}\\
  \AA\transp\vec{y}+\vec{s}&=\;\vec{c}\\
  \|\vec{s}+\mu(z)g(\vec{x})\|_{-\mu(z)g(\vec{x})}
   &\leq\;\b. \end{array}
$$
\end{definition}

The main computational step of each interior-point iteration is to
solve a linearization of the central path equations~\eqref{eq:central}
at the current iterate $z\in \NN_\b$.  The linearization that we will
rely on is as follows.  We will see (cf. Proposition~\ref{geomprop}(f)
below) that for all $\vx,\vs\in \int(\KK)$
there exists a unique {\em scaling point} $w\in \KK$\label{glossary:w} such that
\[
H(w) \vx = \vs.
\]
Given $z = (\vx,\vy,\vs) \in \NN_\b$, the {\em Nesterov-Todd
direction} $(\Delta\vx,\Delta\vy,\Delta\vs)$ is the solution to the
following linearization of \eqref{eq:central}:
\begin{equation}\label{Newton}
\begin{array}{rcl}
\AA \Delta \vx &=& 0 \\
\AA\transp \Delta \vy + \Delta \vs &=& 0\\
\Delta \vx + H(w)^{-1} \Delta \vs &=& -(\mu g(\vs) + \vx).
\end{array}
\end{equation}

By~\cite[Proposition 4.6]{VeraEtAl2007}, the initial point in
step (i) of Algorithm IP below is in the central neighborhood
$\NN_\b$. Furthermore, by \cite[Propositions 4.4 and
4.5]{VeraEtAl2007} if the original pair (P)--(D) is well-posed
(i.e. $C(A)<\infty$), then a point $z\in \NN_\b$ with $\mu(z)$
small enough yields a strict solution to either (P) or (D),
whichever is feasible.  Indeed, by \cite[Theorem
3.1]{VeraEtAl2007} Algorithm IP halts in at most $\Oh(\sqrt{r}
(\log r + \log C(A)))$ iterations and yields a solution to
either (P) or (D).  (See \cite{VeraEtAl2007} for details.)

\begin{quote}
{\bf Algorithm IP}$(A)$

\ch{Let $\b,\d$ be the following constants \label{glossary:d}
$$
\b = \frac{1}{15},\;\; \delta = \frac{1}{45}.
$$}
\begin{description}
\item[(i)]  Let
\[\alpha:= \frac{1}{\sqrt{\ch{2}\bfr}}; \;  M := \frac{\alpha\|Ae\|}{\b};\label{glossary:a}
\]\label{glossary:M:1}
and
\[
\begin{array}{lc}
\vx = (\alpha e, 1, \alpha e, 2M, -\alpha A e )\\
\vy = (0,\frac {M} {\alpha} e,\frac {M} {\alpha^2})\\
\vs = (\frac {M} {\alpha} e, \frac {M}
{\alpha^2}, - \frac {M} {\alpha} e,1,0).
\end{array}
\]
\item[(ii)] If $A\transp y \prec_{K} 0$.
then {\tt HALT} and \\
{\tt return} $y$ as a strictly feasible solution for (D).
\item[(iv)] If \label{glossary:sigma}
$ \sigma_{\min} (A\bar H(x)^{-\um} A^T) > \bfr \mu(z),$ then {\tt HALT} and \\
{\tt return}  $ x + \bar H(x)^{-1}A^T(A\bar H(x)^{-1}A^T)^{-1}x''$
as a strictly feasible solution for (P).
\item[(v)] Set $\bar \mu  := \left(1-\frac{\delta}{\sqrt{\ch{2}\bfr}}\right)
\mu(z).$
\item[(vi)] Compute $\Delta z := (\Delta\vx,\Delta\vy,\Delta\vs)$ by
solving \eqref{Newton} for $\mu = \bar \mu$ and update $z$ by setting
\[
z^+ := z + \Delta z \label{glossary:zplus}
\]
\item[(vii)] Go to (ii).
\end{description}
\end{quote}
\bigskip

It should be noted that the analysis in \cite{VeraEtAl2007} assumes
that all computations are performed with infinite precision.  Our
initial step for a finite-precision algorithm is to show that the
results in \cite{VeraEtAl2007} can be extended to make room for
computational errors.  In particular, Lemma~\ref{lem:keepCentral}
below shows that even if the system~\eqref{Newton} is solved
inexactly, we can still ensure that the iterates remain in the central
neighborhood.

\begin{lemma}\label{lem:keepCentral} Let $z \in \NN_\b$,
$\bar \mu = (1-\frac{\d'}{\sqrt{\bfr}})\mu(z)$ with $|\d'-\d|\leq
\frac{\d}{24}$ and $z^+=z+\D z$ be such that
\begin{equation}\label{eq:keepCentral}
\begin{array}{rl}
 \AA\D\vec{x}&=0\\
 \AA\transp\D\vec{y}+\D\vs&=0\\
 \D\vx+H(w)^{-1}\D\vs&=-(\bar{\mu} g(\vs)+\vx)+\ch{\varrho},
 \end{array}
\end{equation}
where $\|\ch{\varrho}\|\leq \frac{\mu(z)}{120\bfr (2\bfr \mu(z)+1)}$\label{glossary:varrho}. Then
$z^+\in \NN_\b$ and $|\mu(z^+)-\bar \mu| \le \frac{\mu}{120\bfr^2}$.
\end{lemma}
\proof
See Section~\ref{s:PfLKC}.
\eproof

We note that when the solution $\Delta z$ to \eqref{Newton} is computed exactly, \ch{i.e., when $\varrho = 0$ in \eqref{eq:keepCentral},} the point $z^+:= z + \Delta z$ satisfies $\mu(z^+) = \bar \mu$. For details, see~\cite{Renegar99}.

The following two lemmas are in the same spirit as \cite[Propositions
4.4 and 4.5]{VeraEtAl2007}.  In particular, they guarantee that if
either (P) or (D) is strictly feasible then a point $z\in \NN_\b$ with
$\mu(z)$ small enough yields a strict solution to either (P) or
(D). Lemma~\ref{lem:boundSigmaPart} provides the relevant bound for
$\mu(z)$ in the case (P) is strictly feasible and
Lemma~\ref{lem:s_geq_rho} does so for a strictly feasible (D).

\begin{lemma}\label{lem:boundSigmaPart}
Let $z=(\vx,\vy,\vs)\in \NN_\b$ and
assume $\rho_P(A)>0$. Then
\begin{equation}\label{eq:boundSigmaPart}
\sigma_{\min}(A\bar H(x)^{-1}A\transp)  \geq \left(\frac{(1-\b)\rho_P(A)}
  {\b+2\bfr}\right)^2-(2\bfr \mu(z))^2.
\end{equation}
The latter in turn implies that if $\mu(z) <
\frac{(1-\b)\rho_P(A)}{2\bfr(\bfr + \b)}$ then the point
\[
x + \bar H(x)^{-1} A\transp (A \bar H(x)^{-1} A\transp)^{-1} x''
\]
is a strict solution to {\em (P)}.
\end{lemma}

\proof  This is an immediate consequence of the proof of Proposition~4.4
in~\cite[pages 259--260]{VeraEtAl2007}.
\eproof

\begin{lemma}\label{lem:s_geq_rho}
Let $z=(\vx,\vy,\vs)\in \NN_\b$ and assume $\rho_D(A) >0$.
Then $$\|H(s)\| \le \frac{4\bfr^2}{(1-\b)^2 \rho_D(A)^2}.$$
In particular, for $i=1,\dots, r$
\begin{equation}\label{eq:s_geq_rho}
s_{i0}-\|\bar s_i\|\geq \frac{1-\b}{2\bfr\sqrt{r}}\rho_D(A).
\end{equation}


Furthermore, if $\mu(z) < 4\bfr^2(1-\b)\rho_D(A)$ then $y$ is
a strict solution to {\em (D).}
\end{lemma}

\proof
This is an immediate consequence of the proof of Proposition~4.5
in~\cite[page 260]{VeraEtAl2007}.
\eproof

The rest of this section is devoted to proving
Lemma~\ref{lem:keepCentral} and a technical lemma related to the conditioning of the matrix $ $ arising at each interior-point iteration of Algorithm IP. In \S\ref{s:useful} we state and prove a few technical
statements which will be employed in the subsequent proofs. Then
in \S\ref{s:PfLKC} we prove Lemma~\ref{lem:keepCentral}. The proof
of Lemma~\ref{lem:keepCentral} is a straightforward adaptation of
the proof of Theorem 3.7.3 in \cite{Renegar99}. Section~\ref{s:PfLboundSP} presents Lemma~\ref{thm.condnumbers}, which is similar in spirit to
Lemma~\ref{lem:boundSigmaPart}.  This technical result will be crucial in our finite precision analysis in Section~\ref{s:fp}.

\subsection{A few useful relations}\label{s:useful}

The analysis of IPMs heavily relies on the properties of the
barrier function. Here we briefly recall a few essentials that
will be used later. More details can be found in \cite{Renegar99}.
The barrier function $f$ gives rise, for each point $x$ in the
domain $D_f$ of $f$, to a {\em local inner product} $\<\ ,\ \>_x$
induced by $x$ and defined by
$$
  \<u,v\>_x=\<u,H(x)v\>.
$$
The local norm $\|\ \|_x$ is then given by
$\|v\|_x=\<v,v\>_x^{\um}$. In the local inner product $\<\ ,\
\>_x$, the gradient at $y$ is $g_x(y):=H(x)^{-1}g(y)$ and the
Hessian is $H_x(y):=H(x)^{-1}H(y)$.

Our function $f$ defined by \eqref{eq:f} is a self-scaled barrier
with the barrier parameter $\nu= 2\bfr$. We will also use single
components of $f$: $f_i = -\ln(x_{0i}^2 -\|\bar{x_i}\|^2)$ for all
$i = 1, \dots, \bfr$. For each $f_i$ the barrier parameter is $\nu
= 2$. Our development relies on the following key properties of
self-scaled barrier functions \cite{Renegar99}.

\begin{proposition}\label{geomprop}
Let $f$ be a $\nu$-self-scaled barrier function and $x\in D_f\subseteq \R^{\bn}$.
\begin{description}
\item[(a)] If $\|y-x\|_x<1$ then $y\in D_f$ and, for
 all $v \not = 0$,
\[
  1-\|y-x\|_x \leq \frac{\|v\|_y}{\|v\|_x} \leq \frac{1}{1-\|y-x\|_x}.
\]
\item[(b)]
$\{z \in D_f: \ip{z-x}{g(x)} \geq 0\} \subseteq \{z: \|z-x\|_x
\leq \nu\}.$
\item[(c)]
$-g(x) \in D_f$, $-g(-g(x))=x$, $H(-g(x)) = H(x)^{-1}$, and
$\|H(x)^{-1}\| \leq \|x\|^2$.

\item[(d)] For $t > 0$
\[
  g(t x) = \frac{1}{t} g(x),\;\text{ and }
  H(t x) = \frac{1}{t^2}H(x).
\]
\item[(e)] $H(x)x=-g(x)$ and $\ip{x}{g(x)} = - \nu$.
\item[(f)]
Given another point $s\in D_f$, there exist a unique ``scaling
point'' $w\in D_f$ such that
\[
  H(w) x = s, \text{ and } H(w)g(s) = g(x)
\]
and a unique ``reverse scaling point'' $w^*\in D_f$ such that
\[
  H(w^*) s = x, \text{ and } H(w^*)g(x) = g(s).
\]
Furthermore, $w^*:=-g(w)$ and, \ch{for all $\mu > 0$}, the points $\bar
w:=\sqrt{\mu}\,w$ and $\bar w^* := \sqrt{\mu}\,w^*$ satisfy
$$
    \|\bar{w^*}-s\|_{\bar{w^*}}=\|x-\bar{w}\|_{\bar{w}},
$$
and
\[
  \|s + \mu g(x)\|_{-\mu g(x)} \geq \min\left\{\frac{1}{5},
  \frac{4}{5} \|x - \bar w\|_{\bar w}\right\}.
\]
\item[(g)]
If $\|y-x\|_x\leq1$ then $\|g_x(y)-g_x(x)-H_x(x)(y-x)\|_{x}\leq
\frac{\|y-x\|_{x}^2}{1-\|y-x\|_x}$.
\item[(h)] If $\|x-y\|_{y}\leq1$ then
$\|v\|_{-g_y(x)}\leq(1+\|x-y\|_y)\|v\|_y$ for all $v\in\R^{\bn}$.
\eproof
\end{description}
\end{proposition}

\begin{lemma}\label{lem:1}
Let $z\in\NN_\b$ and denote $\mu=\mu(z)$. Then
 $\|x\|=\|x'\|\leq 1$, $\|x''\|\leq\t \leq 2\bfr\mu$,
$\|s''\|\leq1$, $\|s'\|\leq t_s\leq2\bfr\mu$, and $\|s\|\leq
2\bfr\mu+1$.
\end{lemma}
\proof
The bounds on $x$, $x'$ and $s''$ follow from the equalities
$\AA{\vec{x}}=\vec{b}$ and $\AA\transp\vec{y}+{\vec{s}}=\vec{c}$
together with $\vec{x},\vec{s}\in\KK$. The inequalities
$\|x''\|\leq\t\leq2\bfr\mu$ and $\|s'\|\leq t_s=\eta\leq2\bfr\mu$
follow from the equalities
$\t+\eta=\vec{c}\transp\vec{x}-\vec{b}\transp\vec{y}
=\vec{x}\transp\vec{s}=2\bfr\mu$. Since $z\in\NN_\b$, we have
$s=As''-s'$ and therefore
$\|s\|\leq\|s'\|+\|A\|\|s''\|\leq2\bfr\mu+1$. \eproof

The next lemma bounds the norm of the scaling matrix using the results above.

\begin{lemma}\label{lem:normH}
Assume $z \in \NN_\b$, with $\b \le 1/15$. Then
\[
\ch{ \|H(w)\|, \|H(w)^{-1}\| \leq   \frac{4(2\bfr \mu(z)+1)^2}{\mu(z)}.}
\]
\end{lemma}

\proof \ch{Let $\mu = \mu(z)$, $\bar w = \sqrt{\mu} w$, and $\bar w^* = -\sqrt{\mu}g(w).$} 
By Proposition~\ref{geomprop}(f) we have $\|\vec{x} - \overline
w\|_{\overline w} \leq \frac{5}{4}\b$, therefore, by Proposition~\ref{geomprop}(a) and (c), respectively,
\[
 \|H(\overline w)^{-1}\| \leq \frac{1}{(1-\frac{5}{4}\b)^2}
 \|H(\vec{x})^{-1}\| \leq \frac{1}{(1-\frac{5}{4}\b)^2}
\|\vec{x}\|^2.
\]
From Lemma~\ref{lem:1} we have $\|\vx\|^2\leq \ch{3+8\mu^2\bfr^2}$. Hence,
\[
\|H(w)^{-1}\| = \frac{1}{\mu} \|H(\overline w)^{-1}\| \leq
\frac{3+8\mu^2\bfr^2}{(1-\frac{5}{4}\b)^2\mu}\leq
\ch{\frac{4+16\mu^2\bfr^2}{\mu} \le \frac{4(2\bfr \mu+1)^2}{\mu}}.
\]
Similarly, applying Proposition~\ref{geomprop}(f) and Lemma~\ref{lem:1} to $H(\overline w^*)^{-1}$ we have
\begin{equation}\tag*{\qed}
\|H(w)\| = \frac{1}{\mu} \|H(\overline w^*)^{-1}\| \leq
\frac{\|\vs\|^2}{(1-\frac{5}{4}\b)\mu}\leq \frac{4(2\bfr\mu+1)^2}{\mu}.
\end{equation}

\begin{lemma}\label{lem:normg}
\ch{Assume} $z\in \NN_\b$ with $\b\le\frac{1}{15},$ \ch{and $\frac{|\bar \mu - \mu(z)|}{\mu(z)} \le \frac{1}{5\sqrt{2\bfr}}$}. Then
$$
\|H(w)^{-\um}(\bar \mu g(\vec x)+\vec s)\|\leq \frac{\mu(z)^{\um}}{2}.
$$
\end{lemma}

\proof \ch{Let $\mu = \mu(z)$, $\bar w = \sqrt{\mu} w$, and $\bar w^* = -\sqrt{\mu}g(w).$} Since $z\in\NN_\b$,  by
Proposition~\ref{geomprop}(f) we have
\begin{equation}\label{eq:4.62}
\b \geq \|\mu g(\vx)+\vs\|_{-\mu g(\vx)}= \frac{4}{5}\|\vx-\bar
w\|_{\bar w} = \frac{4}{5}\|\bar w^*-\vs\|_{\bar w^*}.
\end{equation}
Then, applying Proposition~\ref{geomprop} and~\eqref{eq:4.62}
twice and using $\b\le1/15$, we obtain
\begin{equation}\label{eq:4.63}
\|\bar \mu g(\vx)+\vs\|_{\bar w^*}\leq \frac{\|\bar \mu
g(\vx)+\vs\|_{\vs}}{1-\frac{5}{4}\b}\leq \frac{\|\bar\mu
g(\vx)+\vs\|_{-\mu g(\vx)}}{(1-\frac{5}{4}\b)(1-\b)}\leq
\frac{5}{3}\|\bar\mu g(\vx)+\vs\|_{-\mu g(\vx)}.
\end{equation}
By the triangle inequality and \eqref{eq:4.62}
\begin{align}\label{eq:4.64}
\|\bar\mu g(\vx)+\vs\|_{-\mu g(\vx)}
 & \leq \frac{|\bar \mu -
\mu|}{\mu}\|g(\vx)\|_{-g(\vx)}+\|\mu g(\vx)+\vs\|_{-\mu
g(\vx)}\leq \frac{4}{15}.
\end{align}
From \eqref{eq:4.63} and \eqref{eq:4.64} we have
\begin{equation}\tag*{\qed}
\|H(w)^{-\um}(\bar \mu g(\vec x)+\vec s)\| = \mu^\um \|\bar \mu
g(\vx)+\vs\|_{\bar w^*} \leq   \frac{4}{9}\mu^\um \le \frac{\mu^\um}{2}.
\end{equation}

\begin{lemma}\label{lem:9} Let $z\in \NN_\b$ for some
$\b\le\frac{1}{15}$. Then for $i=1,\dots, \bfr$
\begin{equation}\label{eq:4.38}
x_i\transp s_i \geq 2(1-\b) \mu(z)
\end{equation}
and
\begin{equation}\label{eq:4.39}
(x_{i0}^2-\|\bar{x_i}\|^2)(s_{i0}^2-\|\bar{s_i}\|^2)\geq 4 (1-\b)^2\mu(z)^2.
\end{equation}
\end{lemma}

\proof \ch{Let $\mu = \mu(z)$.} 
From $\|s_i+\mu g(x_i)\|^2_{-\mu g(x_i)}\leq \b^2$ and
Proposition~\ref{geomprop} we have
\begin{equation}\label{eq:4:40}
\b^2 \geq \|s_i\|^2_{-\mu g(x_i)}-\frac{2}{\mu}\<x_i,s_i\>+2.
\end{equation}
From Proposition~\ref{geomprop}(a,e)
\begin{equation}\label{eq:4.41}
   \|s_i\|_{-\mu g(x_i)}\geq \|s_i\|_{s_i}
   \left(1- \|s_i+\mu g(x_i)\|^2_{-\mu g(x_i)}\right)\geq \sqrt{2}(1-\b).
\end{equation}
Now \eqref{eq:4.38} follows from \eqref{eq:4:40} and \eqref{eq:4.41}.
By Proposition~\ref{geomprop}
$$
   0\leq \|s_i +\mu g(x_i)\|^2_{s_i}= \|s_i\|^2_{s_i}-2\mu \<g(x),g(s)\>
  +\|-\mu g(x_i)\|_{s_i}^2
$$
and by the definition of $f$
$$
\<g(x), g(s)\>= \frac{4\mu\<x_i,s_i\>}{(x_{i0}^2-\|\bar{x_i}\|^2)(s_{i0}^2-\|\bar{s_i}\|^2)}.
$$
Therefore,
\begin{equation}\label{eq:4.42}
(x_{i0}^2-\|\bar{x_i}\|^2)(s_{i0}^2-\|\bar{s_i}\|^2)\geq \frac{8\mu \<x_i,s_i\>}{\|s_i\|^2_{s_i}+\mu^2\|g(x_i)\|_{s_i}^2}.
\end{equation}
By Proposition~\ref{geomprop}
\begin{equation}\label{eq:4.43}
\|-\mu g(x_i)\|_{s_i}\leq \frac{\|-\mu g(x_i)\|_{-\mu g(x_i)}}{1-\|s_i+\mu g(x_i)\|_{-\mu g(x_i)}}\leq \frac{\sqrt{2}}{1-\b}.
\end{equation}
Now \eqref{eq:4.39} follows from \eqref{eq:4.38}, \eqref{eq:4.42} and \eqref{eq:4.43}.
\eproof

\begin{lemma}\label{lem:10} Let $z\in \NN_\b$ \ch{for some
$\b\le\frac{1}{15}$}. Then
\begin{equation}\label{eq:4.44}
\t \geq (1-\b)\mu(z).
\end{equation}
\end{lemma}
\proof By Lemma~\ref{lem:9}, taking $i=r+1$ in \eqref{eq:4.38}, and using $\|s''\|\leq 1$ from Lemma~\ref{lem:1}
$$
2(1-\b) \mu \leq \t\t_s +{x''}\transp s'' \leq \t(1+\|s''\|)\leq 2\t,
$$
which yields \eqref{eq:4.44}.
\eproof

\subsection{Proof of Lemma~\ref{lem:keepCentral}}\label{s:PfLKC}

\ch{Throughout this proof, let $w$ be the scaling point of the pair $\vx$, $\vs$ and  $\bar w = \sqrt{\bar \mu} w$.}
By Lemma~\ref{lem:normH} and the assumptions on the norm of $\varrho$ and on $\bar \mu$ we
have the following bound
\begin{equation}\label{eq:norm.r}
\|\varrho\|_{\bar w} \leq \frac{1}{\bar \mu^\um}\|H(w)^\um\|\|\varrho\|\leq \frac{1}{\bar \mu^\um}\frac{2(2\bfr \mu+1)}{\mu^\um}\frac{\mu}{120(2\bfr \mu+1)}\leq \frac{\mu^\um}{60\bar \mu^\um} \leq\frac{1}{50}.
\end{equation}
Since  $z \in \NN_\b$, and  $\d' \leq
\frac{1}{45}\cdot\frac{25}{24}<\frac{\sqrt{6}}{105}<
\frac{\sqrt{2\bfr}}{105}$, we have
\begin{eqnarray*}
\|\vs + \bar \mu g(\vx)\|_{-g( \vx)}&\leq&
\|\vs + \mu
g(\vx)\|_{-g(\vx)}+|\mu-\bar \mu|\|g(\vx)\|_{-g( \vx)}\\
&\leq& \ch{\b \mu + \delta' \mu <\frac{\bar \mu}{13}. }
\end{eqnarray*}
Therefore, by Proposition~\ref{geomprop}(f)
\begin{equation}\label{eq:8*}
 \| \vx - \bar w\|_{\bar w}\leq \frac{5}{4}\|\vs
 +\bar \mu g(\vx)\|_{-\bar \mu g( \vx)}\leq
 \frac{5}{4}\cdot \frac{1}{13} = \frac{5}{52}.
\end{equation}
Recall that by Proposition~\ref{geomprop} (e)
\begin{equation}\label{eq:0001}
\bar w = -H(\bar w)^{-1}g(\bar w) = -g_{\bar w}(\bar w),
\end{equation}
and that in the local inner product the Hessian $H_{\bar w}(\bar w)$ is the identity
\begin{equation}\label{eq:0002}
H_{\bar w}(\bar w) = H^{-1}(\bar w) H(\bar w)=I.
\end{equation}
Define
\begin{align}\label{def.u}
  u &:= g_{\bar w }(\vx)+2\bar w-\vx\notag\\
    & = g_{\bar w}(\vx)-g_{\bar w }(\bar w)+H_{\bar w}(\bar w)\bar w-H_{\bar w}(\bar w)\vx \quad \text{(by \eqref{eq:0001} and \eqref{eq:0002})}\\
    & = g_{\bar w}(\vx)-g_{\bar w }(\bar w)-H_{\bar w}(\bar w)(\vx- \bar w)\notag.
\end{align}

Since
$
 \|\vx - \bar w\|_{\bar w}< \frac{5}{52} <1,
$
Proposition~\ref{geomprop}(g) yields
\begin{equation}\label{norm.u}
\begin{array}{rcl}
  \|u\|_{\bar w}&=&\|g_{\bar w}(\vx)-g_{\bar w }(\bar w)-H_{\bar w}(\bar w)(\vx- \bar w)\|_{\bar w} \\ [2ex]
  &\le& \frac{\|\vx-\bar w\|^2_{\bar w}}{1-\|\vx-\bar w\|_{\bar w}} \\ [2ex]
  &\leq& \frac{\frac{25}{2704}}{1-\frac{5}{52}}\\[2ex]
  &=&\frac{25}{52\cdot 47}.
  \end{array}
  \end{equation}
Observe that
\begin{eqnarray*}
  \D \vx +H(w)^{-1}\D\vs
  &=& - H(w)^{-1}(\vs +\bar \mu g(\vx))+\ch{\varrho}\\
  &=&-\vx -\bar\mu H(w)^{-1}g(\vx)+\ch{\varrho}
      \qquad \text{by Prop.~\ref{geomprop}(f)}\\
  &=& - \vx - H(\bar w)^{\ch{-1}}g( \vx)+\ch{\varrho}
      \qquad \text{by Prop.~\ref{geomprop}(d)}\\
  &=& -\vx  - g_{\bar w}(\vx)+\ch{\varrho}\\
  &=& 2(\bar w- \vx)-u+\ch{\varrho} \quad(\text{by \eqref{def.u}}).
\end{eqnarray*}
Hence $\bar w - \vx = \frac{1}{2}(\D  \vx
+H(w)^{-1}\D  \vs+u-\ch{\varrho})$ and so
\begin{eqnarray*}
  \bar w - \vx^+ &= &\frac{1}{2}(-\D \vx +H(w)^{-1}\D
 \vs+u-\ch{\varrho}), \quad\mbox{and}\\
  \bar w - H(w)^{-1}\vs^+ &= &\frac{1}{2}(\D\vx
  -H(w)^{-1}\D  \vs+u-\ch{\varrho}).
\end{eqnarray*}
Consequently,
$$
  H(w)^{-1}\vs^+=2\bar w -\vx^+-u+\ch{\varrho}.
$$
Since $\D\vx\perp_w H(w)^{-1}\D\vs$, we have
\begin{eqnarray*}
  \|\bar w - \vx^+\|_{\bar w}&\leq &\frac{1}{2}
  \left(\|-\D\vx +H(w)^{-1}\D
  \vs\|_{\bar w}+\|u\|_{\bar w}+\|\ch{\varrho}\|_{\bar w}\right)\\
  &=&\frac{1}{2}(\|\D \vx + H(w)^{-1}
  \D\vs\|_{\bar w}+\|u\|_{\bar w}+\|\ch{\varrho}\|_{\bar w})\\
  &=&\|\bar w - \vx -\frac{1}{2}u+\frac{1}{2}\ch{\varrho}\|_{\bar
  w}+\frac{1}{2}\|u\|_{\bar w}+\frac{1}{2}\|\ch{\varrho}\|_{\bar w}\\
  &\le&\|\bar w - \vx\|_{\bar w}+\|u\|_{\bar
  w}+\|\ch{\varrho}\|_{\bar w}\\
  &=&\frac{5}{52}+\frac{25}{52\cdot
  47}+\frac{1}{47}=\frac{6}{47}
  \quad(\text{by \eqref{eq:norm.r}, \eqref{eq:8*}, and \eqref{norm.u}}).
\end{eqnarray*}
Thus
\begin{eqnarray*}
\|H(w)^{-1}\vs^++g_{\bar w}(\vx^+ )\|_{\bar w}&=&
\|g_{\bar
w}(\vx^+ )+2\bar w -\vx^+-u+\ch{\varrho}\|_{\bar w}\\
&=& \|g_{\bar w}(\vx^+)+2\bar w -\vx^+\|_{\bar
w}+\|u\|_{\bar
w}+\|\ch{\varrho}\|_{\bar w}\\
&=& \|g_{\bar w}(\vx^+)-g_{\bar w}(\bar w )-H_{\bar w}(\bar
w)^{-1}(\vx^+-\bar w)\|_{\bar w}+\|u\|_{\bar
w}+\|\ch{\varrho}\|_{\bar w}\\
 &\leq & \frac{\|\bar w -
\vx^+\|^2}{1-\|\bar w - \vx^+\|_{\bar w}}\|+\|u\|_{\bar
w}+\|\ch{\varrho}\|_{\bar w}\quad(\text{by Prop.~\ref{geomprop}(g)})\\
&\leq & \frac{36}{47\cdot 41}+\frac{25}{2444}+\frac{1}{50}
<\frac{107}{2132}
\quad(\text{by \eqref{eq:norm.r} and \eqref{norm.u}}).
\end{eqnarray*}
From Proposition~\ref{geomprop}(h) we have for all $v\in\R^{\bn}$
$$
 \|v\|_{-g_{\bar w}(\vx^+)}\leq (1+\|\bar w - \vx^+\|_{\bar
 w})\|v\|_{\bar w}\leq \left(1+\frac{6}{47}\right)\|v\|_{\bar
 w}=\frac{53}{47}\|v\|_{\bar w}.
$$
Thus
\begin{equation}\label{eqn.step1.0}
\begin{array}{rcl}
  \|\vs^++\bar \mu g(\vx^+)\|_{-g(\vx^+)}
   & = &\bar\mu\|H(w)^{-1}\vs^++g_{\bar w}(\vx^+)\|_{-g_{\bar w}
       (\vx^+)} \\
  &< &\frac{53}{47}\cdot
  \frac{107}{2132}\,\bar\mu<\left(\frac{1}{15}-0.01\right)\bar \mu.
\end{array}
\end{equation}
To finish we need to show that \ch{$z^+ = (\vx,\vy,\vs) \in \NN_\b$ and} $\mu(z^+)$ is close to $\bar \mu$. \ch{Since $\|\bar w - \vx^+\|_{\bar w} \le \frac{6}{47} < 1$, Proposition~\ref{geomprop}(g) yields $\vx^+ \in \int(\KK)$ and so $-g(\vx^+) \in \int(\KK).$ From \eqref{eqn.step1.0} we get $\|\vs^++\bar \mu g(\vx^+)\|_{-\bar \mu g(\vx^+)} < 1$ and thus Proposition~\ref{geomprop}(g) again yields $\vs^+ \in \int(\KK)$.
Furthermore,} by assumption we have
$$
\D\vx+H(w)^{-1}\D\vs + \vx=-\bar{\mu} g(\vs)+\ch{\varrho}.
$$
Taking inner product with $\vs$ and using
Proposition~\ref{geomprop}(e,f) we get
\begin{equation}\label{eqn.step1.2}
 \<\vs,\D\vec{x}\>+\<\vx,\D\vec{s}\>+\<\vs,\vx\>=-\bar{\mu}
 \<\vs,g(\vec{s})\>
 +\<\vs,\ch{\varrho}\>=2\bfr\bar{\mu}+\<\vs,\ch{\varrho}\>.
\end{equation}
Since $\D \vx \perp \D \vs$, we have
$
\langle \vx+\D \vx,
\vs +\D \vs\rangle
 =\< \vx, \vs\>+\< \vx,\D
\vs\>+\< \D \vx, \vs\>,
$
so by \eqref{eqn.step1.2}
\begin{equation}\label{eqn.step1.1}
\mu(z^+) = \frac{1}{2\bfr}\langle \vx+\D \vx,
\vs +\D\vs\rangle = \frac{1}{2\bfr}\left(\< \vx, \vs\>+\< \vx,\D
\vs\>+\< \D \vx, \vs\>\right) = \bar \mu +\frac{\langle \vs, \ch{\varrho}\rangle}{2\bfr}.
\end{equation}
Using \eqref{eqn.step1.1}, Lemma~\ref{lem:1} and the assumption on $\ch{\varrho}$, we get
\begin{equation}\label{eq:bound.r2}
|\mu(z^+)-\bar \mu| = \frac{|\<\vs,\ch{\varrho}\>|}{2\bfr}
  \leq \frac{\mu}{120\bfr^2}
  =    \frac{\bar \mu}{120\bfr^2\left(1-\frac{\d'}{\sqrt{2\bfr}}\right)}
  <\frac{\bar \mu}{120}.
\end{equation}
Therefore, from \eqref{eqn.step1.0} and \eqref{eq:bound.r2} we get
\begin{eqnarray}\label{eq:009}
 \|\vs^++\mu(z^+) g(\vx^+)\|_{-g(\vx^+)}&\leq&
 \|\vs^++\bar \mu g(\vx^+)\|_{-g(\vx^+)}+|\mu(
 z^+)-\bar\mu|\|g(\vx^+)\|_{-g(\vx^+)}\notag\\
 &=& \|\vs^++\bar \mu g(\vx^+)\|_{-g(\vx^+)}+2\bfr|\mu(
 z^+)-\bar\mu|\\
 &\leq &\left(\frac{1}{15}-0.01\right)\bar \mu+
 \frac{1}{120}\bar \mu< \frac{1}{15} \mu(z^+).\notag
\end{eqnarray}

\ch{Observe that $z^+$ satisfies the linear equations $\AA\vx^+=0$, $\AA\transp\vy^++\vs^+=0$ by assumption. This together with $\vx,\vs \in \int(\KK)$ and \eqref{eq:009} yields $z^+\in \NN_\b$.}
\eproof

\subsection{On the condition of the matrix $\AA H(w)^{-1/2}\AA\transp$}
\label{s:PfLboundSP}

The purpose of this section is to present Lemma~\ref{thm.condnumbers}
which will be crucial in our finite precision analysis in
Section~\ref{s:fp}.  This lemma is in the same spirit as
Lemma~\ref{lem:boundSigmaPart}. While it will not be used until
Section~\ref{s:fp}, we choose to place it here to facilitate
understanding of the proof.

We will rely on the following key characterization of the distance to ill-posedness due to Renegar~\cite[Theorem 3.5]{Renegar95b}. For a detailed discussion of this and related results, see also~\cite{Lewi99,Lewi01,pena}.
\begin{proposition}\label{prop:infdist}
For any given linear operator $A:\R^n\rightarrow\R^m$ and any cone
$K$
$$
 \rho_P(A)=\sup\{\d:\|v\|\leq\d\Rightarrow v\in\{Ax:\|x\|\leq1,x\in K\}\}
$$
and
\begin{equation}\tag*{\qed}
  \rho_D(A)=\sup\{\d:\|u\|\leq
  \d\Rightarrow u\in\{A\transp y:\|y\|\leq1\}+K^*\}.
\end{equation}
\end{proposition}
We will also rely on the following perturbation result, an extension of~\cite[Theorem 5.1]{VeraEtAl2007}.

\begin{lemma}\label{lem:technicalCN} Let $\b\leq\frac{1}{15}$,
$z = (\vec x, \vec y, \vec s) \in \NN_\b$.  Assume $\bar b\in \R^{\bm}$ is such
that $\langle\vec y, \bar b\rangle\leq 0$ and
\begin{equation}\label{Lem.cond}
\AA\vec v = \bar b \Rightarrow  \left\|\frac{1}{\mu^{\um}}H(w)^\um\vec v\right\|\geq 1
\end{equation}
where $w$ is the scaling point of the pair $\vec x, \vec s$.

If
$\a>\frac{\b+2\bfr}{1-\frac{9}{4}\b}$ then the optimal value of the perturbed problem
\begin{equation}\label{eq:PrtSys}
\begin{array}{rl}
\min    & \vec c^T \vec u \\
s.t.    & \AA \vec u = \vec b +\alpha \bar b,\\
        & \vec u\succeq_{\KK} 0,
\end{array}
\end{equation}
is at least $\vec c^T\vec x$.
\end{lemma}
\proof
This follows by putting together \cite[Theorem 5.1]{VeraEtAl2007} and Proposition~\ref{geomprop}(a,f) as we next explain.
Since $z \in \NN_\b$ Proposition~\ref{geomprop}(f) yields
\[
\|\vx-\bar w\|_{\bar w} \le \frac{5}{4}\b < 1,
\]
where $\bar w = \sqrt{\mu}w$.  Applying Proposition~\ref{geomprop}(a) twice we obtain
\[
\|\vec v\|_{-\mu g(\vs)} \ge (1-\b) \|\vec v\|_{\vx} \ge (1-\b) \left(1-\frac{5}{4}\b\right)\|\vec v\|_{\bar w} \ge \left(1-\frac{9}{4}\b\right)\|\vec v\|_{\bar w}
\]
for all $\vec v$.  Furthermore, observe that $\|\vec v\|_{\bar w} = \left\|\frac{1}{\mu^{\um}}H(w)^\um\vec v\right\|$ for all $\vec v$.  Hence \eqref{Lem.cond} implies that
\[
\AA\vec v = \frac{1}{1-\frac{9}{4}\b}\bar b \Rightarrow  \|\vec v\|_{-\mu g(\vs)}\geq 1.
\]
Therefore, by \cite[Theorem 5.1]{VeraEtAl2007} it follows that the optimal value of \eqref{eq:PrtSys} is at least $\vec c^T\vec x$.
\eproof

\begin{lemma}\label{lem:11}
Let $\D b = (\D b^I,\D b^{II},\D b^{III})\in \R^{m+n+1}=\R^\bm$.
Then there exists $\vec u$ such that
$$
\AA \vec u = \vec b +\D b, \quad \vec u \succeq_\KK 0, \quad
\text{and}\quad \vec c\transp \vec u \ch{\le} \max\{0,
(1+2\rho_P^2(A))^\um\|\D b\|-\rho_P(A)\}.
$$
\end{lemma}
\proof Let \label{glossary:l}
$$
\l = \min\left\{1,\rho_P(A)\frac{1-\|\D b^{II}\|-|\D
b^{III}|}{\|\D b^{I}\|}\right\}.
$$
From Proposition~\ref{prop:infdist} it follows that there exists a
$u$ such that
$$
Au = \l\D b^{I}, \quad u \succeq_\KK 0, \quad \|u\|\leq
\l\frac{\|\D b^I\|}{\rho_P(A)}.
$$
Let $\vec u = (u, 1+\D b^{III}, u+\D b^{II}, (1-\l)\|\D
b^{I}\|,(1-\l)\D b^{I})$. Observe that by construction
$$
\AA\transp u = \vec b -\D b, \quad \vec u \ch{\succeq}_K 0.
$$
Finally,
\begin{align*}
\vec c \transp \vec u
  & = (1-\l)\|\D b^{I}\|\\
  & = \max \{0, \|\D b^I\|+\rho_P(A)(\|\D b^{II}\|+|\D
  b^{III}|-1)\}\\
  & \leq \max \{0, (1+2\rho_P^2(A))^{\um}\|\D b\|-\rho_P(A)\},
\end{align*}
where the last inequality can be obtained by elementary analysis.
\eproof

\begin{lemma}\label{thm.condnumbers}
Let $z\in\NN_\b$. Then if $\rho_P(A)>0$,
\begin{equation}\label{eq:boundSigmaFull}
\sigma_{\min}(\mu^{\um}H(w)^{-1/2}\AA\transp)  \geq
\frac{\mu(z)}{6 \bfr};
\end{equation}
\begin{equation}\label{eq:boundNormB}
\sigma_{\max}(\mu^{\um}H(w)^{-1/2}\AA\transp) \leq
\mu^\um\|H(w)^{-1/2}\|\|\AA\transp\|\leq \ch{4(1 +2\bfr \mu(z))}.
\end{equation}
\end{lemma}
\proof Lemma~\ref{lem:technicalCN} implies that there exists $\D b
\in \R^{\bm} $ satisfying
\begin{equation}\label{eq:4.48}
\|\D b \|\leq
\frac{\b+2\bfr}{1-\frac{9}{4}\b}\sigma_{\min}\left(\mu^\um
H(w)^{-\um}\AA\transp\right)
\end{equation}
and such that the optimal value of the following problem
\begin{equation}\label{eq:PrtSysCpy}
\begin{array}{rl}
\min    & \vec c^T \vec u \\
s.t.    & \AA \vec u = \vec b +\Delta b,\\
        & \vec u\succeq_{\KK} 0
\end{array}
\end{equation}
is at least $\vec c \transp \vx$. Assume that $\|\D b
\|<(1+2\rho_P^2(A))^{-\um}(\vec c \transp \vx +\rho_P(A))$. Then
by Lemma~\ref{lem:11} the optimal value of \eqref{eq:PrtSysCpy} is
$0 <\vec c\transp \vx$, which contradicts the earlier conclusion.
Therefore, $\|\D b\|$ must satisfy
\begin{equation}\label{eq:4.50}
\|\D b \|\geq (1+2\rho^2_P(A))^{-\um}(\vec c \transp
\vx+\rho_P(A)).
\end{equation}
Putting \eqref{eq:4.48} and \eqref{eq:4.50} together, we get
$$
\sigma_{\min}\left(\mu^\um H(w)^{-\um}\AA\transp\right)\geq
\frac{(1-\frac{5}{3}\b)(\vec c \transp
\vx+\rho_P(A))}{(\b+2\bfr)(1+2\rho^2(A))^\um}.
$$
Since $\<\vec c, \vx\>= \t \geq (1-\b)\mu(z)$ by
Lemma~\ref{lem:10}, $0< \rho_P(A)\leq \|A\|\leq 1$ and
$\b<\frac{1}{5}$, we have \eqref{eq:boundSigmaFull}.
Inequality \eqref{eq:boundNormB} follows from the bound \ch{$\|\AA\|\le 2$}
and Lemma~\ref{lem:normH}:
\begin{equation}\tag*{\qed}
\sigma_{\max}(\mu(z)^\um H(w)^{-\um}\AA^T)\leq
\mu(z)^\um\|H(w)^{-\um}\|\|\AA^T\|  \leq \ch{4(1+2\bfr\mu)}.
\end{equation}

\section{Finite precision analysis}\label{s:fp}

\subsection{Floating-point arithmetic}\label{sec:float}

Here we briefly recall the basics of floating-point
arithmetic which we will use in this paper. A slightly more extensive introduction is in~\cite[\S7]{CP01}. Detailed treatments can be found in books on numerical linear algebra such as~\cite{Higham96}.

We call {\em floating-point numbers} a set $\F\subset\R$\label{glossary:F} containing
$0$, {\em rounding map} a transformation $\rou:\R\to\F$\label{glossary:rou} and {\em
round-off unit} a constant $u\in\R$\label{glossary:u} satisfying $0<u<1$. We require
for such a triple that the following properties hold:
\begin{description}
\item[(i)]
For any $x\in\F$, $\rou(x)=x$. In particular $\rou(0)=0$.

\item[(ii)]
For any $x\in\R$, $\rou(x)=x(1+\d)$ with $|\d|\leq u$.
\end{description}

We also define on $\F$ arithmetic operations following the scheme
$$
  x\tilde\circ y=\rou(x\circ y)
$$
for any $x,y\in\F$ and $\circ\in\{+,-,\times,/\}$ so that
$$
  \tilde\circ:\F\x\F\to\F.
$$
It follows from (ii) above that, for any $x,y\in\F$ we have
$$
        x\tilde\circ y=(x\circ y)(1+\d),\qquad |\d|\leq u.
$$
We will also use a floating-point version $\tilde{\sqrt{\ }}$ of the
square root which, similarly, satisfies
$$
        \tilde{\sqrt{x}}=\sqrt{x}(1+\d),\qquad |\d|\leq u.
$$
When combining many operations in floating-point arithmetic,
quantities such as $\prod_{i=1}^n(1+\d_i)^{\rho_i}$ naturally
appear.  The proof of the following propositions can be found in
Chapter~3 of~\cite{Higham96}. The notation they introduce, the
quantities $\g_n$\label{glossary:g} and $\th_n$\label{glossary:th}, and the relations showed therein,
will be widely used in our round-off analysis.

\begin{proposition}\label{propHigham}
If $|\d_i|\leq u$, $\rho_i\in\{-1,1\}$ and $nu<1$ then
$$
  \prod_{i=1}^n(1+\d_i)^{\rho_i}=1+\theta_n
$$
where
$$
        |\theta_n|\leq \g_n=\frac{nu}{1-nu}.
$$
\eproof
\end{proposition}

\begin{proposition}\label{propHigham2}
For any positive integer $k$ such that $ku<1$ let $\theta_k$ be any
quantity satisfying
$$
        |\theta_k|\leq \g_k=\frac{ku}{1-ku}.
$$
The following relations hold.
\begin{description}
\item[1)]
$(1+\theta_k)(1+\theta_j)=1+\theta_{k+j}$,

\item[2)]
$$
\frac{1+\theta_k}{1+\theta_j}=\left\{
\begin{array}{ll}
        1+\theta_{k+j}& \mbox{if $j\leq k$}\\
        1+\theta_{k+2j}& \mbox{if $j> k$,}
\end{array}\right .
$$

\item[3)]
If $ku,ju\leq1/2$ then $\g_k\g_j\leq \g_{\min\{k,j\}}$,

\item[4)]
$i\g_k\leq \g_{ik}$,

\item[5)]
$\g_k+u\leq \g_{k+1}$,

\item[6)]
$\g_k+\g_j+\g_k\g_j\leq \g_{k+j}$.\eproof
\end{description}
\end{proposition}

When computing an arithmetic expression $q$ with a round-off
algorithm, errors will accumulate and we will obtain another
quantity which, we recall, we denote by $\fl(q)$\label{glossary:fl}. We will also write
$\Error(q)=|q-\fl(q)|$\label{glossary:Error}.

An example of round-off analysis which will be useful in the sequel
is given in the next proposition whose proof can be found in
Section~3.1 of~\cite{Higham96}.

\begin{proposition}\label{scalar}
There is a round-off algorithm which, with input $x,y\in\R^n$,
computes the dot product of $x$ and $y$. The computed value
$\fl(\langle x,y\rangle)$ satisfies
$$
   \fl(\langle x,y\rangle)=\langle x,y\rangle+
   \th_{\lceil\log_2 n\rceil+1}\langle |x|,|y|\rangle
$$
where $|x|=(|x_1|,\ldots,|x_n|)$. In particular, if $x=y$ the
algorithm computes $\fl(\|x\|^2)$ satisfying
\begin{equation}\tag*{\qed}
   \fl(\|x\|^2)=\|x\|^2(1+\th_{\lceil\log_2 n\rceil+1}).
\end{equation}
\end{proposition}

The following result deals with summation errors. The proof can be
found in \cite{Higham96}, Section~4.2.

\begin{proposition}\label{sum}
There is a round-off algorithm which, with input $x\in\R^n$,
computes the sum of $x_i$. The computed value
$\fl(\sum_{i=1}^nx_i)$ satisfies
\begin{equation}\tag*{\qed}
   \fl\left(\sum_{i=1}^nx_i\right)=\sum_{i=1}^nx_i+
   \th_{\lceil\log_2 n\rceil}\sum_{i=1}^n|x_i|.
\end{equation}
\end{proposition}

In the next section we will have to deal with square roots. The
following result will help us to do so.

\begin{proposition}\label{sqrtLemma}
Let $\th\in\R$ such that $|\th|\leq 1/2$. Then,
$\sqrt{1+\th}=1+\th'$ with $|\th'|\leq |\th|$. In particular, for
$a\geq 0$
\begin{equation}\label{eq:root}
  \fl\left(\sqrt{a(1+\th_k)}\right)=\sqrt{a}(1+\th_{k+1}).
\end{equation}
\end{proposition}

\proof By the intermediate value theorem we have that
$\sqrt{1+\th}-1=|\th|(\sqrt{\xi})'$ with $\xi\in(1-|\th|,1)$ if
$\th<0$, $\xi\in(1,1+\th)$ otherwise. But
$$
  |(\sqrt{\xi})'|=\left|\frac{1}{2\sqrt{\xi}}\right|\leq \frac{1}{\sqrt2}
$$
the last since $|\xi|\geq 1/2$.

Then \eqref{eq:root} follows from the above. \eproof
Our choice of $u = \phi(\mu(w))$\label{glossary:phi}, for the function $\phi$ in \eqref{eq:u} below, guarantees that
$ku<1/2$ holds whenever we encounter $\th_k$, and consequently,
$\th_k \leq 2ku.$ We will therefore not bother the reader by
repeating this fact each time we use it.
\medskip

\subsection{The finite precision algorithm}

In this section we present a finite precision algorithm that determines which one
of (P) or (D) is strictly feasible and provides a solution. In the case when the dual  problem (D) is feasible, after sufficiently refining the precision we will be able to obtain an
exact feasible solution to (D), however, for the primal problem
only an approximation to a feasible solution is possible due to
the structure of the problem: we cannot compute a point on the
linear subspace $Ax=0$ exactly with finite precision. However, we
can obtain a forward-approximate primal solution of any desired
accuracy. To describe this in more detail, we need the following
definition of a $\gamma$-approximate solution.

\begin{definition}
Let $\gamma\in(0,1)$. A point $\hat{x}\in\R^n$ is a {\em
$\gamma$-forward solution} of the system $Ax=0$, $x\succ_{K}0$, if
$\hat{x}\succ_{K}0$, and there exists $\breve{x}\in\R^n$ such that
$$
  A\breve{x}=0,\qquad
  \breve{x}\succ_{K}0
$$
and
$$
  \|\hat{x}-\breve{x}\|\leq\gamma\|\hat{x}\|.
$$
The point $\breve{{x}}$ is said to be an {\em associated solution}
for $\hat{x}$. A point is a {\em forward-approximate solution} of
$Ax=0$, $x\succeq_{K}0$, if it is a $\gamma$-forward solution of
the system for some $\gamma\in(0,1)$.  \ch{Observe that by definition, the existence of a $\gamma$-forward solution automatically guarantees the existence of a strict solution.} 
\end{definition}

We are now ready to present our main result and give a precise description of the related algorithm. The proof of Theorem~\ref{main_th} is deferred to Section~\ref{s:main_proof}.

\begin{theorem}\label{main_th}
There exists a finite precision algorithm which, with input a matrix
$A\in\R^{m\x n}$ and a number $\gamma\in(0,1)$, finds either a
strict $\gamma$-forward solution $x\in\R^n$ of $Ax=0$,
$x\succeq_{K}0$, or a strict solution $y\in\R^m$ of the system
$A\transp y\preceq_{K}0$. The machine precision varies during the
execution of the algorithm. If (P) is strictly feasible, the
finest required precision is
$$
   u^*=\left(\bfc(n+m)^{5/2}r^{8} C(A)^{7/2}\left(1+\frac{1}{\gamma}\right)^{7/2}\right)^{-1},
$$
and in the case when (D) is strictly feasible,
$$
   u^*=\left(\bfc(n+m)^{5/2}r^{11.5} C(A)^{7/2}\right)^{-1},
$$
where $\bfc$ is a universal constant\label{glossary:bfc:1}. The number of main
(interior-point) iterations of the algorithm is bounded by
$$
  \Oh\left(r^{1/2}(\log(r)+\log(C(A))+|\log\gamma|)\right)
$$
if (P) is strictly feasible and by the same expression without the
$|\log\gamma|$ term if (D) is.
\end{theorem}


\begin{remark}
In the numerical analysis literature, fixed precision is used more
commonly than variable precision. We note here that from our
variable precision analysis we can obtain a fixed
precision one. Indeed, assume the precision $u$ is fixed.
Then our algorithm could run until the point
at which it should get a precision finer than $u$. If it found
the answer before this point it could return it (and this answer
would be guaranteed to be correct).  If not, it could halt and
return a failure message. Furthermore,
the only reason for $u$ to be insufficient is that
$C(A)$ is too large.  Solving the bound for $u$ in
Theorem~\ref{main_th} we
obtain a lower bound $C_u$ for $C(A)$. Thus, the failure
message could be something like
{\em ``The condition of the data is larger than $C_u$. To solve
the problem I need more precision.''}  \ch{We note that although the statement of Theorem~\ref{main_th} depends on the condition number $C(A)$, Algorithm FP described below does not require any information on $C(A)$ as input.  The only required input are the matrix $A$ and a constant $\gamma \in (0,1)$.} 
\end{remark}

We are now ready to describe our primal-dual algorithm. This is
essentially an extension of Algorithm IP from
Section~\ref{s:ipm} with some additional features.  One of these
features is the stopping criteria and the other one is the presence
of finite precision and the adjustment of this precision as the
algorithm progresses. To ensure the correctness of the algorithm,
the precision will be set to
\begin{equation}\label{eq:u}
 \phi(\mu(z)) :=  \frac{\mu(z)^{7/2}}{\bfc\bfr\bn^{5/2}(2\bfr\mu+1)^{11/2}}
\end{equation}
at each iteration. Here $\bfc$ is a universal constant.

Let $\b=\frac{1}{15}$ and $\d = \frac{1}{45}$.
\begin{quote}
\begin{description}
 \item{\ } {\bf Algorithm FP} $(A,\gamma)$
 \item[(i)] Set the machine precision to
$u := \frac{1}{\bfc \bfr^7(\bm+\bn)^{5/2}}$\\
$\a :=\frac{1}{\sqrt{2\bfr}},M=\frac{\a\|Ae\|}{\b}$\label{glossary:M:1FP}\\
$z:=\left(\a e,1,\a e,2M,-\a
Ae,0,\frac{M}{\a}e,-\frac{M}{\a^2},\frac{M}{\a}e,\frac{M}{\a^2},-\frac{M}{\a}e,1,0\right)$
 \item[(ii)] Set the machine precision to $u := \phi(\mu(z))$.
 \item[(iii)] If for $i=1,\ldots,r$\\
$s_{i0}-\|\bar s_i\|-6\mu(z)\bfr >0$\\
then {\tt HALT} and {\tt return} $y$ as a strict solution for
$A\transp y \preceq_{K} 0$.
 \item[(iv)] If
$ \sigma_{\min} (\bar H(x)^{-\um} A^T)\ge \frac{3\bfr \mu(z)}{\g}, $
then {\tt HALT} and \\
{\tt return} $x$ as a $\gamma$-forward solution for $Ax=0,\;
x\succeq_{K}0$.
 \item[(v)] Set $\bar \mu  :=
\left(1-\frac{\d}{\sqrt{2\bfr}}\right) \mu(z).$
 \item[(vi)] Update
$z$ by solving the linearization~\eqref{Newton} of \refer{eq:central} for $\mu =
\bar \mu$. \item[(vii)] Go to (ii).
\end{description}
\end{quote}
\bigskip

The matrix $\bar H(w)$ used in step (iv) is the upper-left
$n\times n$ block of $H(w)$, where $w$ is the scaling point
of $(\vec{x},\vec{s})$.

The precise way we solve the system in (vi) is as follows:
\begin{description}
\item[(a)] Compute a solution $\D\vec{y}$ of
\begin{equation}\label{red.eqns}
  (\AA H(w)^{-1}\AA\transp)\Delta\vec{y}=\AA H(w)^{-1}
 (\vec{s}+\bar{\mu} g(\vec{x})).
\end{equation}

\item[(b)] Let $\vec{y}:=\vec{y}+\D\vec{y}$,
$$
  \left(\begin{array}{l}\D x\\
\D\t\end{array}\right):=\left(\begin{array}{lllll}\Id_n&0&0&0&0\\
0&0&0&1&0\end{array}\right)(
H(w)^{-1}\AA\transp\Delta\vec{y}-(\bar{\mu} g(\vec{s})+\vec{x})).
$$
Then set $x:=x+\D x$ and $\t:=\t+\D\t$.

\item[(c)] Let
$$
\begin{array}{l}x':=x\\
 t:=1\\
 x'':=-Ax\\
 s:=y'-A\transp y\\
\end{array}
\qquad
\begin{array}{l}
 s':=-y'\\
 t_s:=\eta\\
 s'':=-y\\
 \t_s=1.\end{array}
$$
\end{description}

\begin{remark}\label{rem:DisregrdAddition}
The finite-precision errors in the computations in
(b) and (c) are negligible compared to the errors involved in
solving the linear system on step (a).  Therefore, for ease of
exposition, we will assume that the computations in (b) and (c)
in step (vi) are exact.  We also assume that the initial point $z$ in step (i) and the value of $\bar \mu$ in step (v) of Algorithm FP are computed exactly.  We stress that these assumptions have no
consequences in the complexity or accuracy bounds.
By making them we can greatly reduce the length
of our exposition and focus our analysis on the
critical stages of the algorithm. \ch{We assume that the smallest singular value in step (iv) above is computed using a backward stable algorithm (e.g., QR
factorization). This guarantees that the computed $\fl(\sigma_{\min}(\bar H(x)^{-\um} A^T))$ is the
exact $\sigma_{\min}((\bar H(x)^{-\um} A^T)+E)$ for a matrix $E$ with $ \|E\|\leq
\constCQR n^{5/2}u\|\bar H(x)^{-\um} A^T\|$\label{glossary:constCQR} for some universal constant
$\constCQR$ (see, e.g., \cite[Chapter 2]{Bjor96}).}

\end{remark}

Under the assumption of infinite precision on steps (b) and (c) the next point $z^+$ thus defined lies in the linear subspace
\{$\AA\vx=\vec{b}$, $\AA\transp \vy -\vs = \vec{c}$\}.  Moreover, $\D z = (\D x, \D y, \D z)$ satisfies system
\eqref{eq:keepCentral} for some (possibly large) $r$.

The crux of our finite precision analysis is the estimation of the floating-point errors in step (a) above, which we present in Section~\ref{sec.finite.Newton}.   That analysis relies on the following technical lemma.

\begin{lemma}\label{B}
Assume $z\in\NN_\b$ and let $w= w(z)$ be its scaling point. With
precision $u = \phi(\mu(z))$ we can compute $B =
H(w)^{-\um}\AA\transp$ and $D = \bar H(w)^{-\um}A\transp$ (where $\bar
H(w)$ is the upper-left $n\times n$ block of $H(w)$) satisfying
\begin{equation}\label{eq:Bound_B_D}
   \|\fl(B)-B\|,\|\fl(D)-D\|\leq \frac{1}{336 \cdot 240}\cdot
   \frac{\mu^2}{\bfr(2\bfr \mu +1)^2}
\end{equation}
as well as $q = -H(w)^{-\um}(\bar \mu g(\vec x)+\vec s)$
satisfying
\begin{equation}\label{eq:Bound_Err_q}
\|\fl(q) - q\|\leq \frac{1}{16\cdot 240}\cdot
\frac{\mu^{3/2}}{(2\bfr \mu +1)^2}.
\end{equation}
\end{lemma}

The proof of Lemma~\ref{B} in turn relies on the following technical result.

\begin{lemma}\label{lem:err_det} Let $z\in \NN_{\b}$, and
the finite-precision computations are performed with
$u=\phi(\mu (z))$. Then
\begin{equation}\label{eq:err_det}
\Error\left(\frac{\det s_i}{\det x_i}\right)\leq \frac{\det s_i}{\det x_i}\g_M, \qquad  \Error\left({\det s_i}{\det x_i}\right)\leq {\det s_i}{\det x_i}\g_M,
\end{equation}
where\label{glossary:det}
$$
M = \frac{(2\bfr\mu+1)^2(\log_2(m+n)+2)}{\mu^2(1-\b)^2}.\label{glossary:M:2}
$$
\end{lemma}
\proof Observe that from Proposition~\ref{scalar}
\begin{equation}\label{eq:011}
\Error(\det s_i) = \|s_i\|^2\th_{\lceil\log_2 n_i\rceil+1},\quad \Error(\det x_i) = \|x_i\|^2\th_{\lceil\log_2 n_i\rceil+1}
\end{equation}
for all $i\in \{1,\dots, \bfr\}$. Let $\kappa_i: = \lceil\log_2
n_i\rceil+1$. From \eqref{eq:011} we have
\begin{align*}
\Error\left(\frac{\det s_i}{\det x_i}\right)
  & \leq  \frac{\det s_i+\|s_i\|^2\g_{\kappa_i}}{\det x_i-\|x_i\|^2\g_{\kappa_i}}(1+\g_1)
    -\frac{\det s_i}{\det x_i}\\
  & =  \frac{\det s_i}{\det x_i}\cdot\frac{\g_1\det s_i \det x_i
    +\|s_i\|^2\g_{\kappa_i+1}\det x_i+\|x_i\|^2\g_{\kappa_i}\det s_i }{(\det x_i-\|x_i\|^2\g_{\kappa_i})
     \det s_i}\\
  &\leq  2\frac{\det s_i}{\det x_i}\left(\frac{\|x_1\|^2}{\det x_i}\g_{\kappa_i}
     +\frac{\|s_i\|^2}{\det s_i}\g_{{\kappa_i}+1}+\g_1\right),
\end{align*}
and from Lemma~\ref{lem:9}
$$
\det x_i \geq \frac{4\mu^2(1-\b)^2}{\|s_i\|^2}, \qquad \det s_i \geq \frac{4\mu^2(1-\b)^2}{\|x_i\|^2}.
$$
Therefore, using Lemma~\ref{lem:1}
\begin{align*}
  \Error\left(\frac{\det s_i}{\det x_i}\right)
   & \leq
  \frac{\det s_i}{\det x_i}\cdot \frac{\|x_i\|^2\|s_i\|^2\g_{\kappa_i}
   +\|x_i\|^2\|s_i\|^2\g_{\kappa_i+1}+\g_1}{2\mu^2(1-\b)^2}\\
   & \leq
  \frac{\det s_i}{\det x_i}\cdot \frac{(2\bfr \mu+1)^2\g_{\kappa_i}
   +(2\bfr \mu+1)^2\g_{\kappa_i+1}+\g_1}{2\mu^2(1-\b)^2}\\
   & \leq
  \frac{\det s_i}{\det x_i}\cdot \frac{(2\bfr \mu+1)^2}{\mu^2(1-\b)^2}\g_{\kappa_i+1},
\end{align*}
which yields the first inequality in \eqref{eq:err_det}. The
second relation is obtained analogously. \eproof

\proofof{Lemma~\ref{B}}
It is well-known (see \cite[\S3.2]{Tsuchiya}) that the scaling matrix
$H$ has a block-diagonal structure, where each block corresponds
to a Lorentz cone; moreover, each individual block can be represented
as follows 
$$
   H_i(w(x_i,s_i))^{-1/2}=\gscmx^{-1}
   \left(\begin{array}{cc}\a&-\z\transp\\
     -\z&I+\frac{\z\z\transp}{1+\a}\end{array}\right),\label{glossary:z}
$$
where
$\gscmx=\left[\frac{s_{i0}^2-\|\bar{s_i}\|^2}{x_{i0}^2
-\|\bar{x_i}\|^2}\right]^{1/4}$, $\a=\frac{\xi_0}{\det(\xi)^{\um}}$,
and $\z=\frac{\bar{\xi}}{\det(\xi)^{\um}}$ with
$\xi=(\xi_0,\bar{\xi})=(\gscmx^{-1}s_{i0}+\gscmx
x_{i0},\gscmx^{-1}\bar{s_i}-\gscmx\bar{x_i})$\label{glossary:xi} for all $i=1,\dots, \bfr$.

From Proposition~\ref{scalar} for all $i=1,\dots,\bfr$, $
  \Error(\<x_i,s_i\>)\leq \|x_i\|\|s_i\|\g_{\log_2(n+m)+1}$,
and using Lemmas~\ref{lem:1} and~\ref{lem:9}
\begin{equation}\label{eq:0003}
\Error(\<x_i,s_i\>)\leq \<x_i,s_i\>\frac{2\bfr \mu+1}{2(1-\b)}\g_{\log_2(n+m)+1}\leq \<x_i,s_i\>\gamma_{M} .
\end{equation}
By Proposition~\ref{sqrtLemma} and
Lemma~\ref{lem:err_det}
\begin{equation}\label{eq:0004}
\fl (\sqrt{\det x_i \det s_i}) = \sqrt{\det x_i \det s_i} (1+\gamma_{M+1}).
\end{equation}
Therefore, from \eqref{eq:0003} and \eqref{eq:0004}
\begin{align*}
\Error (\det \xi)
& = \Error \left(2\left(\sqrt{\det x_i \det s_i}+\<s_i,x_i\>\right) \right)\\
& = \left|2\left(\sqrt{\det x_i \det s_i} (1+\theta_{M+1})+\<x_i,s_i\>(1+\theta_{M})\right)(1+\theta_2)-\det \xi\right|\\
& = \left|2\left(\sqrt{\det x_i \det s_i} \theta_{M+1}+\<x_i,s_i\>\theta_{M}\right)(1+\theta_2)+\det \xi \theta_2\right|\\
& = \det \xi \left|\theta_{M+1}+\theta_2+\theta_{M+1}\theta_2\right| \leq \det\xi \gamma_{M+3},
\end{align*}
the last inequality due to Proposition~\ref{propHigham2} 4).

Further
\begin{equation}\label{eq:B10}
\Error(\gscmx^{-1}\a) = \Error\left(\frac{\gscmx^{-2}s_{i0}+x_{i0}}{\sqrt{\det
\xi}}\right)= \gscmx^{-1}\a\th_{3M+4};
\end{equation}
\begin{equation}\label{eq:B11}
\fl(\gscmx^{-1}\z)_j = \fl\left(\frac{\gscmx^{-2}s_{ij}-x_{ij}}{\sqrt{\det
\xi}}\right)=
(\gscmx^{-1}\z)_j+\frac{\gscmx^{-2}|s_{ij}|+|x_{ij}|}{\sqrt{\det
\xi}}\th_{3M+4},
\end{equation}
hence
$$
\Error ((\gscmx ^{-1}\z)_j)\leq \gscmx^{-1}\a \g_{3M+4}.
$$
It remains to evaluate the errors in the bottom-left block of
$H_i(w(x_i, s_i))$. We have $\Error (\xi_k)\leq \xi_0\g_{2M+4}$,
then
$$
\Error \left(\frac{\z_k\z_l}{1+\a}\right) = \Error
\left(\frac{\xi_k\xi_l}{\det \xi^{\um}+\xi_0}\right)\leq
\frac{\xi^2_0}{\det \xi^\um +\xi_0}\g_{5M+15}
=\frac{\z^2_0}{1+\a}\g_{5M+15}
$$
and
$$
\Error \left(\gscmx^{-1}\frac{\z_k\z_l}{1+\a}\right)\leq
\gscmx^{-1}\frac{\z^2_0}{1+\a}\g_{11M+33};
$$
$$
\Error \left(\gscmx^{-1}+\gscmx^{-1}\frac{\z_k\z_l}{1+\a}\right)\leq
\gscmx^{-1}\left(1+\frac{\z^2_0}{1+\a}\right)\g_{11M+33}.
$$
Finally, we have
$$
\|\fl (H(w)^{-\um}) -H(w)^{-\um} \|\leq
(n_i+1)\gscmx^{-1}\max\left\{\a \g_{3M+4},
\left(1+\frac{\z^2_0}{1+\a}\right)\g_{11M+34}\right\}.
$$
Observe that by Lemma~\ref{lem:9}
$$
\det \xi = 2(\sqrt{\det x_i \det s_i}+\<x_i, s_i\>)\geq 8
(1-\b)\mu;
$$
$$
\xi_0 = \gscmx^{-1}s_0 +\gscmx x_0 \leq
\frac{x_0^2+s_0^2}{\sqrt{2\mu(1-\b)}}; \qquad \gscmx^{-1} \leq
\frac{s_0}{\sqrt{2\mu(1-\b)}}.
$$
Then
\begin{equation}\label{eq:alpha}
 \a = \frac{\xi_0}{\det\xi^{\um}}\leq
 \frac{s_0^2+x_0^2}{4\mu(1-\b)};
\end{equation}
$$
1+\frac{\z_0^2}{1+\a} = 1+\frac{\z_0^2}{\det \z+\z_0}\leq
1+\z_0\leq 1+\frac{x_0^2+s_0^2}{4\mu(1-\b)}.
$$
Observe that $x_0^2+s_0^2 <2(2\bfr \mu+1)^2$. Hence we have
$$
 \gscmx^{-1}\max\left\{\a \g_{3M+4},
 \left(1+\frac{\z^2_0}{1+\a}\right)\g_{11M+34}\right\}\leq
 \frac{5(2\bfr \mu +1)^{7/2}}{\mu^{3/2}}\g_{11M+34}.
$$
Therefore,
\begin{equation}\label{eq:0005}
  \|\fl (H(w)^{-\um}) -H(w)^{-\um} \|\leq \frac{5(n+m)(2\bfr \mu
  +1)^{7/2}}{\mu^{3/2}}\g_{11M+34}
\end{equation}
and
$$
  \|\fl (\bar H(w)^{-\um}) -\bar H(w)^{-\um} \|\leq \frac{5n(2\bfr
  \mu +1)^{7/2}}{\mu^{3/2}}\g_{11M+34}.
$$
Now we estimate the error in computing $B = H(w)^{-\um}\AA\transp$.
First, observe that the error for multiplication of
$\fl(H(w)^{-\um})$ by $\AA\transp$ can be estimated as follows (see
\cite[Chapter 22]{Higham96})
$$
  \|\fl ( H(w)^{-\um})\AA\transp -\fl\left(\fl(
  H(w)^{-\um})\AA\transp\right)\|\leq \bn^{2} u \|\fl (
  H(w)^{-\um})\|\|\AA\transp\|+\constCL u^2\label{glossary:constCL}
$$
\ch{for some universal constant $\constCL$.}
Therefore, recalling that $\|{\cal A}\|\leq \ch{2}$ (see
\eqref{eq:normAA}),
\begin{align*}
\|B - \fl(B) \| &
  \leq \| H(w)^{-\um}\AA\transp - \fl(H(w)^{-\um})\AA\transp\|\\
  & \qquad + \|\fl(H(w)^{-\um})\AA\transp-\fl\left(\fl(H(w)^{-\um})\AA\transp\right)\|\\
  & \leq \left(\| H(w)^{-\um} - \fl(H(w)^{-\um})\| + \bn^{2} u\|
  \fl(H(w)^{-\um})\|\right)\|\AA\transp\|+\constCL u^2\\
  & \leq \ch{2}\| H(w)^{-\um} - \fl(H(w)^{-\um})\|\\
  & \qquad + \ch{2}\bn^{2} u\left[ \| H(w)^{-\um} - \fl(H(w)^{-\um})\| +\|(H(w)^{-\um})\|\right]+\constCL u^2\\
  & =\ch{2}\left(1+\bn^{2} u\right)\| H(w)^{-\um} - \fl(H(w)^{-\um})\|
  +\bn^{2} u\|(H(w)^{-\um})\|+\constCL u^2\\
  & < \frac{\ch{40}(n+m)^2(2\bfr
  \mu+1)^{7/2}}{\mu^{3/2}}\g_{12M+34}+\constCL u^2,
\end{align*}
where the last bound follows from \eqref{eq:0005} and
Lemma~\ref{lem:normH}. Since our precision $u$ satisfies
\eqref{eq:u}, this yields \eqref{eq:Bound_B_D} for $B$.  Observe
that here we only care about the order of the problem-driven
parameters, not the constants, as $\bfc$ in the precision update
formula \eqref{eq:u} can be adjusted to accommodate any
multiplicative constants.  The corresponding bound for $D$ is
obtained analogously.

It remains to evaluate the errors in computing $q$. By
straightforward computation we obtain for each `block'
$$
  q_{i0}= -\frac{1}{2}\det \xi_i^\um\left(1-\bar \mu
  \frac{\gscmx_i^{-2}}{\det x_i}\right);
$$
$$
  \bar q_{i}= -\frac{\left(1-\bar \mu \frac{\gscmx_i^{-2}}{\det
  x_i}\right)}{\det\xi_i^\um+\xi_{i0}}\left[\left(s_{i0}+\frac{1}{2}\gscmx_i\det
  \xi_i^\um\right)\bar x_i+\left(\frac{1}{2}\gscmx_i^{-1}\det
  \xi_i^\um+x_{i0}\right)\bar s_i\right].
$$
We obtain \eqref{eq:Bound_Err_q} by a similar argument as when
evaluating $\|\fl(H(w)^{-\um})- H(w)^{-\um}\|$. We omit this
tedious exercise for the sake of brevity. \eproof

%
%

\subsection{Finite-precision analysis of solving the Newton system}\label{sec.finite.Newton}

The main result of this section is Lemma~\ref{thm.reduced},
which bounds the round-off error in the computation in the
reduced equations~(\ref{red.eqns}) in (a), when it is performed
with finite precision.

 Observe that the reduced system of
equations~(\ref{red.eqns}) is equivalent to the least-squares
problem
\[
\min_v \; \|Bv + q \|^2
\]
for $B = H(w)^{-1/2}\AA\transp$ and $q = -H(w)^{-\um}(\bar \mu
g(\vec x)+\vec s)$. We rely on this equivalence to obtain the
bound in Lemma~\ref{thm.reduced}. More precisely, we apply a
known round-off error result for the least-squares problem,
namely Proposition~\ref{prop:golub}. Lemma~\ref{thm.reduced}
follows from Proposition~\ref{prop:golub} and suitable bounds
on the norm of $q$ and singular values of $B$ obtained earlier
in Lemmas~\ref{lem:normg} and \ref{thm.condnumbers}
respectively.

Recall the following stability property of Golub's method for
linear least-squares (cf.~\cite[Chapter~16]{LawHan74}).

\begin{proposition}\label{prop:golub}
Let $\BB\in\R^{p\times l}$ $(p\geq l)$\label{glossary:BB} have full rank. Let $u$
denote the machine precision. If Golub's method is applied to
$$
  \min _{v\in\R^l}\|\BB v + f\|
$$
the computed solution is the exact solution to a problem
$$
  \min_{v\in\R^l}\|(\BB+\d\BB)v+(f+\d f)\|
$$
where
$$
  \|\d\BB\|\leq \constCGolub upl^{3/2}\|\BB\|,\qquad
  \|\d f\|\leq \constCGolub upl\|f\|\label{glossary:constCGolub}
$$
and $\constCGolub $ is a universal constant independent of $p$
and $l$. \eproof
\end{proposition}

\begin{lemma}\label{thm.reduced}
Let $z\in\NN_\b$. With precision $u=\phi(\mu)$ in all
arithmetic operations, we can compute a vector $\fl(\Delta
\vec{y})$ such that
$$
\|(\AA H(w)^{-1} \AA\transp) \fl(\Delta \vec{y}) - \AA
H(w)^{-1}(\bar{\mu} g(\vec{x})+\vec{s})) \| \leq
\frac{\mu(z)}{\ch{120}\bfr (2\bfr\mu(z)+1)}.
$$
In addition, $\fl(\D \vy)\leq 12 \bfr\mu^{-\um} $.
\end{lemma}
\proof Let $B = H(w)^{-\um}\AA\transp$ and $q = H(w)^{-\um}(\bar
\mu g(\vx) +\vs)$, then the system of equations
$$
(\AA H(w)^{-1}\AA\transp)\D\vy = \AA H(w)^{-1}(\bar \mu g(\vx)
+\vs)
$$
can be written as $B\transp B\D\vy = -B\transp q$ and its
solutions are those of the least squares problem
\begin{equation}\label{eq:4.65}
\min_{v\in \R^\bm}\|Bv+q\|.
\end{equation}
Hence we can apply Golub's method to \eqref{eq:4.65} to compute a
solution to the original equation. Let $\D\vy$ be the vector
actually computed by Golub's method when solving \eqref{eq:4.65}.
Then $\D\vy$ is the exact solution of
\begin{equation}\label{eq:4.66}
\min_{v\in \R^\bm}\|\tilde Bv+\tilde q\|
\end{equation}
for some $\tilde B$ and $\tilde q$ satisfying
\begin{equation}\label{eq:4.67}
\|\tilde B - \fl(B)\|\leq \constCGolub u \bm^{3/2}\bn
\|\fl(B)\|, \quad \|\tilde q- \fl(q)\|\leq \constCGolub u\bm
\bn \|\fl(q)\|,
\end{equation}
where $\constCGolub$ is a universal constant. Let $\D B =\tilde
B - B$, $\D q = \tilde q - q$. Since $\D \vy$ is an exact solution
of the least squares problem \eqref{eq:4.66}, we have $\tilde B
\transp \tilde B \D \vy +\tilde B\transp \tilde q = 0$, and thus
\begin{equation}\label{eq:4.68}
B\transp B\D \vy +B\transp q = - (\tilde B\transp \D B +\D
B\transp B)\D \vy - (\tilde B\transp \D q+\D B\transp q).
\end{equation}
From Lemmas~\ref{lem:normH} and \ref{B} and $\|\AA\transp \|\leq
\ch{2}$ we have
\begin{equation}\label{eq:4.69}
\|\fl(B)\|\leq \|H(w)^{-\um}\|\|\AA\transp\|+\|\fl(B) - B\|\leq
\frac{\ch{4}(2\bfr \mu+1)}{\mu^\um}.
\end{equation}
Analogously, Lemmas~\ref{lem:normg}  and \ref{B} yield
\begin{equation}\label{eq:4.70}
\|\fl(q)\|\leq \|q\|+\|\fl(q) - q\|\leq
\frac{\mu^\um}{2}+\frac{\mu^{3/2}}{(2\bfr \mu+1)^2}.
\end{equation}
Then from \eqref{eq:4.67},  \eqref{eq:4.69}, \eqref{eq:4.70} and
our choice of $u$ we have
\begin{equation}\label{eq:4.71}
\|\tilde B - \fl(B)\|\leq
\ch{4}\constCGolub u\bm^{3/2}\bn\frac{2\bfr\mu+1}{\mu^\um}\leq
\frac{1}{336\cdot \ch{240}}\cdot\frac{\mu^2}{\bfr (2\bfr \mu+1)^2};
\end{equation}
\begin{equation}\label{eq:4.72}
\|\tilde q - \fl(q)\|\leq \ch{4}\constCGolub u\bm\bn\mu^\um\leq
\frac{1}{16\cdot \ch{240}}\cdot\frac{\mu^{3/2}}{\bfr (2\bfr \mu+1)^2}.
\end{equation}
\ch{Here we assume that the constant $\bfc$ in \eqref{eq:u} is chosen so that inequalities \eqref{eq:4.71} and \eqref{eq:4.72} hold.  This ensures that the rest of the proof goes through.}
Applying Lemma~\ref{B} again and using \eqref{eq:4.71} and
\eqref{eq:4.72},
\begin{equation}\label{eq:4.73}
\|\D B\|\leq \|\fl(B) - B\|+\|\tilde B - \fl (B)\| \leq
\frac{1}{168 \cdot \ch{240}}\cdot \frac{\mu^2}{\bfr (2\bfr \mu+1)^{2}};
\end{equation}
\begin{equation}\label{eq:4.74}
\|\D q\|\leq \|\fl(q) - q\|+\|\tilde q - \fl (q)\| \leq \frac{1}{8
\cdot \ch{240}}\cdot \frac{\mu^{3/2}}{(2\bfr \mu+1)^2}.
\end{equation}
Then from \eqref{eq:4.73}, \eqref{eq:4.74} and using the bounds on
$\|B\|$ and $\|q\|$ discussed above,
\begin{equation}\label{eq:4.75}
\|\tilde B \|\leq \|B\| +\|\D B\|\leq 4\frac{2\bfr \mu
+1}{\mu^\um}; \qquad \|\tilde q \|\leq \|q\| +\|\D q\|\leq\mu^\um.
\end{equation}
From \eqref{eq:4.73} and \eqref{eq:4.75} we have
\begin{equation}\label{eq:4.76}
\|\tilde B\transp \D B +\D B\transp B \|\leq (\|\tilde B\| +\|\D
B\|)\|\D B\|\leq \frac{1}{\ch{240}}\cdot\frac{\mu^{3/2}}{12\bfr (2\bfr
\mu+1 )}.
\end{equation}
From \eqref{eq:4.73}, \eqref{eq:4.74} and \eqref{eq:4.75}
\begin{equation}\label{eq:4.77}
\|\tilde B\transp \D q+\D B\transp q\|\leq \|\tilde B \|\|\D q\|
+\|\D B\|\|q\| \leq \frac{1}{\ch{240}}\cdot \frac{\mu}{2\bfr \mu+1}.
\end{equation}
It remains to bound $\|\D \vy\|$.  Using
Lemma~\ref{thm.condnumbers} and \eqref{eq:4.73} we have
$$
\sigma_{\min}(\tilde B)\geq \sigma_{\min}(B) -\|\D B\|\geq
\frac{\mu}{6\bfr}- \frac{\mu}{12\bfr}  = \frac{\mu}{12\bfr}.
$$
From \eqref{eq:4.75} we have $\|\tilde q\|\leq \mu^\um$.  Observe
that since $\| \tilde B \D \vy+\tilde q\| =\min_{v}\|\tilde B v
+\tilde q\|$,
\begin{equation}\label{eq:4.78}
\|\D \vy\|\leq \frac{\|\tilde q\|}{\sigma_{\min }(\tilde B)}\leq
12\mu^{-\um}.
\end{equation}
Finally, we have from \eqref{eq:4.68}, \eqref{eq:4.76},
\eqref{eq:4.77} and \eqref{eq:4.78}
\begin{align}
\|B\transp B \D \vy+B\transp q\| & \leq \|\tilde B \transp \D B
+\D B \transp B\|\cdot \|\D \vy\|+\|\tilde B \transp \D q+\D B
\transp q\|\notag\\
&\leq \frac{1}{\ch{240}}\cdot \frac{\mu^{3/2}}{12\bfr (2 \bfr
\mu+1)}\cdot 12\bfr \mu^{-\um}+\frac{1}{\ch{240}}\cdot
\frac{\mu}{2\bfr\mu +1}\notag\\
&\leq \frac{\mu}{\ch{120}\bfr (2\bfr \mu+1)}.\tag*{\qed}
\end{align}
\smallskip

\subsection{Finite-precision analysis of termination conditions}

\begin{lemma}[Dual termination]\label{prop:p}
Let $\rho_D(A)>0$ and  $z\in \NN_\b$ with
$\mu(z)\leq \frac{\rho_D(A)}{40\bfr^3}$. Then
\begin{equation}\label{eq:013}
  \fl \left(s_{i0} -\|\bar {s_i}\|\right)
  >\fl\left(6\mu(\ubar z)\bfr\right) \; \text{for } \; i\in 1:r.
\end{equation}
Moreover, if $z \in \NN_\b$ satisfies \eqref{eq:013}, then the subcomponent $y$ of $z$ is a strict feasible solution to (D); in other words, $A\transp y \prec_K 0$.
\end{lemma}

\proof
From our choice of precision $u = \phi(\mu(z))$ and the fact that $z \in \NN_\b$ it readily follows that
\begin{equation}\label{error.so}
\Error(s_{i0} -\|\bar {s_i}\|) \le \bfr \mu(z)
\end{equation}
On the other hand, by Lemma~\ref{lem:s_geq_rho} we have
\[
s_{i0} -\|\bar {s_i}\| \ge \frac{1-\b}{2\bfr\sqrt{r}}\rho_D(A) \ge
\frac{1-\b}{2\bfr\sqrt{r}}(40 \bfr^3 \mu(z))  \ge 10 \mu(z) \bfr.
\]
Thus
\[
\fl(s_{i0} -\|\bar {s_i}\|) \ge \fl(8\bfr \mu(z)) > \fl(6\bfr \mu(\ubar z)) .
\]
Now assume \eqref{eq:013} holds.  Again, by \eqref{error.so} we get
\[
s_{i0} -\|\bar {s_i}\| \ge 5\bfr \mu(z).
\]
Since  $\AA \vy + \vs = \vec{c}$, in particular $A\transp y - y' + s = 0$. So $-A\transp y = s - y'$. Since $\|y'\| \le \eta \le \tau + \eta = \vec{c}\transp \vx - \vec{b}\transp \vy = 2\bfr \mu(z)$, it follows that for $i=1:r$ we have
\[
s_{i0} - y'_{i0} -\|\bar {s_i} - \bar {y'_i}\| \ge s_{i0} -\|\bar {s_i}\|  - y'_{i0} - \|\bar {y'_i}\| \ge 5\bfr \mu(z) - 2\|y'\| \ge \bfr \mu(z) > 0.
\]
Therefore, $s - y' \succ_K 0$ and consequently $A\transp y = y'-s \prec_K 0.$
\eproof

\begin{lemma}[Primal termination]\label{prop:primal_stop}
Assume $z\in \NN_\b$. If $\mu(z)\leq
\frac{\rho_P(A)}{10\bfr^2}\left(1+\frac{1}{\g}\right)^{-1}$ then
in step (iv) the algorithm yields
\begin{equation}\label{eq:055}
\fl\left(\sigma_{\min} \left(\bar H(x)^{-\um} A\transp\right)\right)\ge
\fl\left(\frac{3\bfr\mu(z)}{\g}\right).
\end{equation}
Moreover, if $z \in \NN_\b$ satisfies \eqref{eq:055} then the subcomponent $x$ of $z$ is a $\g$-forward solution of $Ax=0$,
$x\succeq_{K}0$, and
$$
  \breve{x} = x-\bar H(x)^{-1}A\transp(A \bar
  H(x)^{-1}A\transp)^{-1}Ax
$$
is an associated solution for $x$.
\end{lemma}
\proof Let $D=\bar H(x)^{-\um}A\transp$ and assume that we compute
$\sigma_{\min}(D)$ using a backward stable algorithm (e.g., QR
factorization). Then the computed $\fl(\sigma_{\min} (D))$ is the
exact $\sigma_{\min}(\fl(D)+E)$ for a matrix $E$ with $ \|E\|\leq
\constCQR n^{5/2}u\|\fl(D)\|$ for some universal constant
$\constCQR$ (see \cite[Chapter 2]{Bjor96}). We have
\begin{align}\label{error.svd}
\Error(\sigma_{\min}(D))
 &\le \|\fl(D) - D\|+\constCQR n^2
  u(\|D\|+\|\fl(D)-D\|)\\
 & \leq  (1+\constCQR n^2 u)\|\fl(D) - D\|+\constCQR n^2 u\|\overline H ^{-\um}(x)\|\|A\|\\
 & \leq \frac{ (1+\constCQR n^2 u)\mu^2(z)}{336\cdot 240\cdot \bfr(2\bfr \mu(z)+1)^2}+\constCQR n^2
 u\frac{2+3\mu(z)\bfr}{\mu(z)^{\um}},
\end{align}
where the last inequality follows from Lemmas~\ref{lem:normH}
and \ref{B}. By our choice of $u$ we have
$$
\Error(\sigma_{\min}(D))\le  \bfr\mu(z).$$
Therefore,
$$
 \fl(\sigma_{\min}(D)) \geq \sigma_{\min}(D) -
  \bfr\mu(z).
$$
Using the bound from Lemma~\ref{lem:boundSigmaPart} we have
\begin{align*}
 \fl(\sigma_{\min}(D)) &  >
 4\bfr\mu(z)\left(\frac{\rho_P(A)}{10\bfr^2\mu(z)}-1\right).
\end{align*}
By our choice of $\mu(z)$ this yields $$\fl(\sigma_{\min}(D))> \frac{4\bfr\mu(z)}{\g} \ge \fl\left(\frac{3\bfr\mu(z)}{\g}\right).$$
Now assume \eqref{eq:055} holds.  Again by \eqref{error.svd} we get
\[
 \sigma_{\min}(D) \ge \fl(\sigma_{\min}(D))  -
\bfr\mu(z) \ge  \frac{2\bfr\mu(z)}{\gamma}.
\]
Denote $\D x =
-\bar H(x)^{-1}A\transp(A \bar H(x)^{-1}A\transp)^{-1}Ax$. From
Lemma~\ref{lem:1} and $\|Ax\|=\|x''\|$ we have
\begin{align*}
  \|\D x\|^2_x
     & =(Ax)\transp(A \bar H(x)^{-1}A\transp)^{-1}Ax
     \leq \frac{\|x''\|^2}{\sigma_{\min}(D)^2}
     \leq \frac{(2\bfr \mu(z))^2}{(2\bfr \mu(z)/ \gamma)^2}
     \leq \g^2.
\end{align*}
Furthermore, by Proposition~\ref{geomprop}(a,c) we have
$$
\frac{\|\Delta x\|}{\| x\|} \leq \frac{\|\Delta x\|}{\|H(
x)^{-\um}\|} \leq \|H( x)^{\um}\Delta x\| = \|\Delta x\|_{ x} \leq \gamma < 1.
$$
Therefore, $\breve{x} = x + \Delta x$ is a $\g$-forward solution of $Ax=0$, $x\succeq_{K}0$.
\eproof

\section{Proof of the main result}\label{s:main_proof}

We are finally in a position to prove our main result
(Theorem~\ref{main_th}). We first prove that on every step the
algorithm keeps up with the central path, and at the same time the
value of $\bar \mu$ decreases by a fixed factor. Then we show that
once $\mu$ is small enough to satisfy either dual or primal
termination conditions (Steps (iii) and (iv) of the
algorithm), the algorithm terminates and yields a correct answer.
Then the bound on the number of iterations follows trivially from
the termination bounds on $\mu$ and the factor of the decrease of
$\mu$. Similarly we obtain bounds for the finest precision based
on the precision update function $\phi$ and the aforementioned
bounds on $\mu$.

Before we go ahead with the proof, we need to
guarantee that the initial point satisfies all the necessary
bounds. The following result is an immediate consequence of~\cite[Proposition 4.6]{VeraEtAl2007}.

\begin{lemma}[Computation of the initial point]\label{lem:init}
The initial point
$$
z:=\left(\a e,1,\a e,2M,-\a
Ae,0,\frac{M}{\a}e,-\frac{M}{\a^2},\frac{M}{\a}e,\frac{M}{\a^2},-\frac{M}{\a}e,1,0\right)
$$
where $\a =\frac{1}{\sqrt{2\bfr}}$ and $M=\frac{\a\|Ae\|}{\b}$, satisfies $z\in \NN_\b$ and $\mu(z) = \frac{\alpha \|Ae\|}{\beta} = O(1)$.

\eproof
\end{lemma}

\proofof{Theorem~\ref{main_th}}
We first disregard the halting steps (iii) and (iv) of the algorithm and prove that, no matter how many iterations we have performed, all our iterates stay close enough to the central path. We use an induction argument.

 The induction base is given by Lemma~\ref{lem:init}. We now assume that at the start of Step (ii) of the algorithm the value of $z$  satisfies $z\in \NN_\b$. We need to show that the vector $z^+$ computed in step (vi)
 is also in $\NN_\b$.

From Lemma~\ref{thm.reduced} it follows that the point $z^+$
computed with finite precision in (a),  and infinite precision in (b) and (c), satisfies
$$
\begin{array}{rl}
 \AA\D\vec{x}&=0\\
 \AA\transp\D\vec{y}+\D\vs&=0\\
 \D\vx+H(w)^{-1}\D\vs&=-(\bar{\mu} g(\vs)+\vx)+\ch{\varrho},
 \end{array}
$$
for some $\ch{\varrho}$ with $\|\ch{\varrho}\|\leq \frac{\mu(z)}{120\bfr (1+\mu(z))}$.
Hence from Lemma~\ref{lem:keepCentral} we have $z^+\in \NN_\b$.
Hence, $z\in \NN_\b$ on every iteration.

Now we show the bounds for the number of iterations and the finest
precision.

From Lemma~\ref{prop:primal_stop} we know that once $\mu(z)$
reaches the lower bound of
$\frac{\rho_P(A)}{10\bfr^2}\left(1+\frac{1}{\g}\right)^{-1}$, the
algorithm yields a correct $\g$-approximate solution to the
primal problem.

It follows from Lemma~\ref{lem:keepCentral}  that
$$
\mu(z^+)\leq \bar \mu+\frac{\mu(z)}{120\bfr^2} = \left(1 - \frac{\delta}{\sqrt{2\bfr}}+\frac{1}{120\bfr^2}\right)\mu(z).
$$
Therefore, taking into account that $\bfr^2\geq \sqrt{2\bfr}\cdot
3\sqrt{\frac{3}{2}} > 3 \sqrt{2\bfr}$ and that $\d = 1/45$, we
have
$$
\mu(z^+)<\left(1-\frac{1}{60\sqrt{2\bfr}}\right)\mu(z).
$$
Given an initial value of $\mu(z_0)$, after $k$ iterations we have
$$
\mu(z_k) \leq \left(1-\frac{1}{60\sqrt{2\bfr}}\right)^k \mu(z_0).
$$
Since we want $\mu(z) \leq
\frac{\rho_P(A)}{10\bfr^2}\left(1+\frac{1}{\g}\right)^{-1}$, we
have the condition
$$
\left(1-\frac{1}{60\sqrt{2\bfr}}\right)^k \mu(z_0) \leq
\frac{\rho_P(A)}{10\bfr^2}\left(1+\frac{1}{\g}\right)^{-1}.
$$
By taking logarithms on both sides, and using $\mu(z_0) = O(1)$, we get the desired relation
$$
k = \Oh\left(r^{1/2}(\log(r)+\log(C(A))+|\log\gamma|)\right).
$$

Since the algorithm halts once $\mu(z)\leq
\frac{\rho_P(A)}{10\bfr^2}\left(1+\frac{1}{\g}\right)^{-1}$, we
deduce from our precision update formula \eqref{eq:u} that the
finest required precision $u^*$ satisfies
$$
u^* \geq \left(\bfc(n+m)^{5/2}\bfr^8
C(A)^{7/2}\left(1+\frac{1}{\gamma}\right)^{7/2}\right)^{-1}.
$$

Similarly, for the dual feasible case from Lemma~\ref{prop:p} we
know that $\mu(z)\leq \frac{\rho_D(A)}{40\bfr^3}$ guarantees
successful termination of the algorithm. Hence we similarly get
the bound
$$
k = \Oh\left(r^{1/2}(\log(r)+\log(C(A)))\right),
$$
and for the finest precision we have
\begin{equation}\tag*{\qed}
u^* \geq \left(\bfc(m+n)^{5/2}\bfr^{11.5} C(A)^{7/2}\right)^{-1}.
\end{equation}

\ch{\section*{List of main symbols used in the paper}}

\ch{\begin{tabular}{ll}
$A$, $K$, $r$& page \pageref{glossary:A}\\
 $\rho_P(A)$, $\rho_D(A)$, $C(A)$& page \pageref{glossary:rhoP} \\
 $\AA$, $\vec b$, $\vec c$, $\KK$, $\vs$, $s'$, $s''$, $t$, $t_s$, $\vx$, $x'$, $x''$, $\vy$, $y'$, $\eta$& page \pageref{glossary:KK}\\
 $\bm$, $\bn$, $\bfr$ & page \pageref{glossary:bm}\\
 $e, f(\vec x)$, $\bar f(x)$, $g(\vec x)$, $\bar g(x)$, $H(\vec x)$, $\bar  H(x)$ & page \pageref{eq:f} \\
 $\NN_\b$, $w$, $z$, $\mu(z)$& page \pageref{glossary:NNb}\\
 $\a$, $\b$, $\d$, $\sigma$& page \pageref{glossary:b} \\
 $\rou$, $u$, $\g_k$, $\th_k$& page \pageref{glossary:rou}\\
 $\Error$, $\fl$& page \pageref{glossary:Error}\\
 $\phi(\mu(z))$& page \pageref{glossary:phi} \\
 $\bfc$& page \pageref{glossary:bfc:1}\\
 $\constCQR$& page \pageref{glossary:constCQR}\\
 $\constCL$& page \pageref{glossary:constCL}\\
 $\constCGolub$& page \pageref{glossary:constCGolub}\\
\end{tabular}
}

\end{document}